\documentclass[12pt,twoside]{amsart}

\usepackage{amssymb}

\newtheorem{theorem}{Theorem}[section]
\newtheorem{lemma}[theorem]{Lemma}
\newtheorem{proposition}[theorem]{Proposition}
\newtheorem{corollary}[theorem]{Corollary}
\newtheorem{conjecture}[theorem]{Conjecture}

\newtheorem{problem}[theorem]{Problem}

\theoremstyle{definition}
\newtheorem{definition}[theorem]{Definition} 
\newtheorem{remark}[theorem]{Remark}
\newtheorem{example}[theorem]{Example}

\newtheorem{notation}[theorem]{Notation}

\newcommand{\fa}{{\mathfrak a}}
\newcommand{\fb}{{\mathfrak b}}
\newcommand{\fc}{{\mathfrak c}}

\newcommand{\Ann}{\operatorname{Ann}}
\newcommand{\Soc}{\operatorname{Soc}}

\newcommand{\mif}{\mbox{if} ~}

\newcommand{\reg}{\operatorname{reg}}

\newcommand {\ZZ}{\mathbb{Z}}

\newcommand {\PP}{\mathbb{P}}

\begin{document}

\title[Minimal Resolution of Relatively Compressed Level Algebras]
{Minimal resolution of Relatively Compressed Level Algebras}

\author[J.\ Migliore, R.\ Mir\'o-Roig, U.\  Nagel]{Juan C.\ Migliore,  Rosa M.\
Mir\'o-Roig, Uwe Nagel}

\address{Department of Mathematics,
           University of Notre Dame,
           Notre Dame, IN 46556,
           USA}

\email{Juan.C.Migliore.1@nd.edu}

\address{Facultat de Matem\`atiques,
     Department d'Algebra i Geometria,
     Gran Via des les Corts Catalanes 585,
     08007 Barcelona,
     SPAIN }

\email{miro@ub.edu }

\address{Department of Mathematics,
University of Kentucky, 715 Patterson Office Tower, Lexington, KY
40506-0027, USA}

\email{uwenagel@ms.uky.edu}

\date{\today}

\begin{abstract}
A relatively compressed algebra with given socle degrees is an
Artinian quotient 
$A$ of a given graded algebra $R/\fc$, whose Hilbert function is maximal
among such quotients with the given socle degrees.  For us $\fc$ is usually a
``general'' complete intersection and we usually require that $A$ be level.
The precise value of the Hilbert function of a relatively compressed algebra
is open, and we show that finding this value is equivalent to the Fr\"oberg
Conjecture.

We then turn to the minimal free resolution of a level algebra relatively
compressed with respect to a general complete intersection.  When the algebra
is Gorenstein of even socle degree we give the precise resolution.  When it
is of odd socle degree we give good bounds on the graded Betti numbers.  We
also relate this case to the Minimal Resolution Conjecture of Musta\c t\v a
for points on a projective variety.

Finding the graded Betti numbers is essentially equivalent to determining to
what extent there can be redundant summands (i.e.\ ``ghost terms'') in the
minimal free resolution, i.e.\ when copies of the same $R(-t)$ can occur in
two consecutive free modules.  This is easy to arrange using Koszul
syzygies; we show that it can also occur in more surprising situations that
are not Koszul.  Using the equivalence to the Fr\"oberg Conjecture, we show
that in a polynomial ring where that conjecture holds (e.g.\ in three
variables), the possible non-Koszul ghost terms are extremely limited.

Finally, we use the connection to the Fr\"oberg Conjecture, as well as the
calculation of the minimal free resolution for relatively compressed
Gorenstein algebras, to find the minimal free resolution of general Artinian
almost complete intersections in many new cases.   This greatly extends
previous work of the first two authors.

\end{abstract}

\thanks{The authors thank the Mathematical Research Institute at
Oberwolfach, where
the authors first discussed this problem.  Part of the work for this
paper was done
while the first author was sponsored by the National Security Agency
under Grant
Number MDA904-03-1-0071.  The second author was partially sponsored by
the grant BFM2001-3584.}

\maketitle
\tableofcontents


\section{Introduction} \label{section-intro}

Let $k$ be a field and let $R = k[x_1,\dots,x_n]$
be the homogeneous polynomial ring.  We will say that an Artinian
$k$-algebra $A = R/I$  has {\em 
  socle degrees} 
$(s_1,\ldots,s_t)$ if the minimal generators of its socle (as
$R$-module) have degrees $s_1 \leq \ldots \leq s_t$. Thus, the number
of $s_j$'s that equal $i$  is the
dimension of the component of the socle of $A$ in degree $i$.  For  fixed
socle degrees, a graded Artinian algebra of maximal Hilbert function among
all graded Artinian algebras with that socle degrees is said to be {\em
compressed}.  We extend this notion as follows:

\begin{definition} \label{def-compr}
Let $\fc \subset R$ be a homogeneous ideal and let $0 \leq s_1 \leq
\ldots \leq s_t$ be integers.
   Then a graded Artinian $k$-algebra $A$ is said to be
   {\it relatively compressed with respect to $\fc$} and with socle
   degrees 
   $(s_1,\ldots,s_t)$ if $A$ has maximal length among all graded Artinian
   $k$-algebras $R/I$ satisfying
\begin{itemize}
\item[(i)] $\Soc R/I \cong \oplus_{i=1}^t k(-s_i)$; \item[(ii)]
$\fc \subset I$.
\end{itemize}
\end{definition} 

Equivalently, $A = R/I$ is relatively compressed with respect to $\fc$
if it is a quotient of $R/\fc$ having maximal length and the
prescribed socle degrees. This is a slight extension (allowing $\fc$ to
be Artinian or 
non-saturated) of the notion of ``relatively compressed algebras"
introduced in \cite{CI}, Definition 2.2.  The paper \cite{Iarrobino-compr}
introduces a notion of an algebra being ``compressed relative to an Artin
algebra," but this is unrelated to our notion.

In almost all of our work, $\fc$
will be a complete intersection  (only in Section \ref{sec-exist} will we
extend this to allow $\fc$ to be Gorenstein).  {\em Note that the complete
intersection itself is not necessarily Artinian.}  In many situations it is
important to look for ideals that  contain a regular sequence (i.e.\ complete
intersection) in certain degrees, and to ask for such ideals that are
maximal in some way.  The lex-plus-powers conjecture is an example of such a
problem (cf.\
\cite{richert}).  We are interested in a similar situation, seeking an ideal
that contains a regular sequence of fixed degrees, has  fixed socle degrees,
and has maximal Hilbert function among all such ideals.

In Section \ref{sec-exist} of this paper we give a good
motivation for studying such ideals by showing the connection to the famous
conjecture of Fr\"oberg on the Hilbert function of an ideal of general forms
of fixed degrees.  There is a natural guess for the Hilbert function of a
relatively compressed Artinian algebra, based on an upper bound coming from
the theory of inverse systems (see for instance \cite{tonynotes},
\cite{Iarrobino-LNM}).  However, we also give examples to show that the
``natural guess'' for the Hilbert function of a relatively compressed
algebra need not hold, and we also show that even the choice of the field
(in positive characteristic) can affect this Hilbert function.  We also show
how the choice of the complete intersection can affect this, again in
positive characteristic.  We do not know if this can happen in
characteristic zero, but anyway for many of our results we will assume that
$\fc$ is a ``general'' complete intersection of fixed generator degrees.

The question of when a given Hilbert function can fail to exist for a
given socle degrees, despite it satisfying the natural guess mentioned
above, has been considered elsewhere.  See for instance \cite{CI},
\cite{GHMS}, \cite{zanello}.  Our focus on algebras relatively
compressed with respect to a complete intersection, and our consequences
of this fact, are new.  This part of our work can be viewed as a partial
answer to the question asked near the end of Remark 2.9 of
\cite{CI}, to determine an upper bound for the Hilbert function of
a relative compressed algebra, and to try to see when it may be sharp.

For most of this paper, however, we are interested in {\em level} graded
Artinian algebras, i.e.\ in the case where all socle degrees are
equal, so the socle is 
concentrated in one degree. The study of level algebras was
initiated by Stanley \cite{stanley2}.   Level graded algebras   play 
an important role
in many parts of commutative algebra, algebraic geometry and algebraic
combinatorics.  For instance, a sufficiently general set of points in
projective space is often level (it depends on the number of points) -- cf.\
\cite{lauze}.  Even the Gorenstein case, which is just a special case of
level algebras, has an extensive literature.   See \cite{GHMS} for an
extensive bibliography and overview of level algebras.

The level graded Artinian algebras of maximal Hilbert function among all
level graded algebras of given codimension and socle degree are called {\em
compressed level algebras} and they fill a non-empty Zariski open set in the
natural parameter space. If the socle has dimension one and occurs in degree
$s$ then the algebra $R/I$ is Gorenstein and $I$ can be identified with the
ideal consisting of  partial derivatives of all orders annihilating a general
polynomial $f \in R$ of degree $s$.  The notion of a relatively compressed
Gorenstein algebra naturally arises by requiring that a priori some
partial derivatives of $f$ vanish.

Beyond finding the Hilbert function, a much more subtle question
is to  understand all of the syzygies of a relatively compressed level
Artinian algebra, i.e.\ to find the minimal free
$R$-resolution of $A$ or, equivalently, the graded Betti numbers.  A central
part of this is to determine if there can be redundant summands (i.e.\
``ghost terms'') in consecutive free modules in the minimal free
resolution.  Sometimes this is easy to force with Koszul syzygies.  The
interesting situation is when there are syzygies that cannot be explained
by Koszul relations, the so-called {\em non-Koszul ghost terms}.

Note that even in the case of compressed level Artinian algebras,
very little is known about the minimal free resolution.  In \cite{boij2}
Corollary 3.10, Boij showed that there is a well-defined notion of
``generic'' Betti numbers for compressed level algebras of fixed socle
degrees, and in  Conjecture 3.13 he guessed what they may be.  The
main point is that there should be no  ghost terms.
The first case is that of Gorenstein algebras.  When the socle degree
is even, the result was well-known, following from the almost purity of the
minimal free resolution.  In the case of odd socle degree, the result was
shown by  the first and second authors in \cite{MMR2} Proposition 3.13, as
long as the  initial degree is sufficiently large.  Very little is known
beyond this.

In Section \ref{sec-Gor}, our first main result is to give the precise
minimal free resolution of a relatively compressed Artinian Gorenstein
algebra, $A$, of even socle degree.  Here we assume that $k$ has
characteristic zero, or else that the characteristic satisfies a certain
numerical condition (see Remark~\ref{characteristic remark}).  We find that
any ghost terms in the minimal free resolution of $A$ occur either directly
because of Koszul relations, or indirectly because of Koszul relations and
duality.  We also give a similar result for relatively compressed Artinian
Gorenstein algebras of odd socle degree, but here we are not able to give
the precise resolution (but we show where there is uncertainty).  However, we
are able to give the precise resolution in odd socle degree when the
embedding dimension is 4, the socle degree is odd, and $A$ is relatively
compressed with respect to a general quadric.  We also give a nice
connection to the Minimal Resolution Conjecture for points on complete
intersection varieties (as special case of a conjecture of Musta\c t\v a),
showing that if the conjecture holds then we can find the minimal free
resolution of a general  Gorenstein Artinian algebra of odd socle degree,
relatively compressed with respect to a general complete intersection of
codimension $\leq n-2$.

Section \ref{sec-Conj} deals with level algebras of socle dimension  $\geq 2$.
Our main goal is to see how it can happen that the minimal free resolution
has non-Koszul ghost terms.  We give some conditions that
force such ghost terms.  We also give several examples and conjectures.
Finally, we show that if $R$ satisfies Fr\"oberg's Conjecture then the
minimal free resolution of a relatively compressed level algebra can have
non-Koszul ghost terms only in a very limited way.  This holds, for example,
if $n=3$.

In Section \ref{sec-Applications} we go in the opposite direction.  It is
known that in characteristic zero (and slightly more generally) an ideal of
$n+1$ general forms satisfies Fr\"oberg's conjecture.  It has been
conjectured that an ideal of $n+1$ general forms has the ``expected" minimal
free resolution in the sense that the Betti numbers are the minimal ones
consistent with the Hilbert function (i.e.\ no ghost terms).  The first and
second authors showed this in several cases in \cite{MMR2} (and also gave
some counterexamples).    Here, using our result for Gorenstein algebras in
Section
\ref{sec-Gor}, we show that an ideal of $n+1$ general forms (with generator
degrees satisfying certain conditions) must have the predicted graded
Betti numbers, extending the known results.

The authors would like to thank Tony Iarrobino for helpful comments
about the exposition of this paper and the connections to his work.


\section{Relatively compressed algebras and Fr\"oberg's conjecture}
\label{sec-exist}

Throughout this paper we will use the following notation:

\begin{notation}
Let $k$ be an infinite field (often making further assumptions, such as
characteristic zero).  Let
$R=k[x_1, \cdots, x_n]$ and let $A$ be a graded $k$-algebra. The Hilbert
function of $A$ is denoted by $h_A(t) := \dim_k A_t.$
\end{notation}

In this section we give some basic results about  relatively
compressed algebras.
The main purpose is to establish the connection between them and Fr\"oberg's
conjecture on the Hilbert function of an ideal of general forms.
However, we also
discuss the failure of the ``expected" Hilbert function to occur even
with generic choices. 

If the Artinian algebra $A$ has socle degrees $(s_1,\ldots,s_t)$ then 
$s_t = \max\{s_1,\ldots,s_t\}$ is called {\em the socle degree} of $A$. It
equals the Castelnuovo-Mumford regularity of $A$. 
Moreover, by a {\em general} homogeneous polynomial of degree $d$ we
  mean a polynomial in a suitable Zariski open and dense subset of
  $R_d$. Similarly, a general complete intersection is generated by
  general polynomials. 

We now begin with a simple remark.

\begin{remark} \label{rem-cone} Let $A = R/I$ be an Artinian algebra
    with  socle degrees 
    $(s_1,\ldots,s_t)$. Let $\fc \subsetneq I$ be an Artinian Gorenstein
    ideal. Denote by $e$ the socle degree of $R/\fc$.
    Let $J := \fc : I$ be the residual ideal. Then there are homogeneous
    forms $G_1,\ldots,G_t \in R$ of degree $e-s_1,\ldots,e-s_t$ such
    that
$$
J = \fc + (G_1,\ldots,G_t).
$$
In fact, this follows from the standard mapping cone procedure that
relates the resolutions of $I, J, \fc$.

We will say that $e-s_1,\ldots,e-s_t$ are the {\it expected degrees}
of the extra generators of $J$, i.e.\ of the minimal generators of $J$
that are not in $\fc$.
\end{remark}

We would like to generalize this remark. Let $A = R/I$ be an Artinian algebra
   with  socle degrees 
   $(s_1,\ldots,s_t)$. Let $\fc \subsetneq I$ be a Gorenstein
   ideal of codimension $c$. Denote by $e$ the Castelnuovo-Mumford
   regularity of $R/\fc$ and put
$$
s := \max \{s_1,\ldots,s_t\}.
$$
Let $F_1,\dots,F_{n-c} \in I$ be general forms of degree
$s+1$. Then
$$
\fc' := \fc + (F_1,\dots,F_{n-c})
$$
is an Artinian Gorenstein ideal with socle degree
$$
e' = e + (n-c) s
$$
because $I_{s+1} = R_{s+1}$.

\begin{remark}
If $\fc$ is already Artinian then $c = n$, $\fc = \fc'$, $e = e'$ and there
are no forms
$F_i$ needed.
\end{remark}

Now we link.
   Let $J := \fc' : I$ be the residual ideal. Then there are homogeneous
   forms $G_1,\ldots,G_t \in R$ of degree $e'-s_1,\ldots,e'-s_t$ such
   that
$$
J = \fc' + (G_1,\ldots,G_t) = \fc + (F_1,\dots,F_{n-c}) + (G_1,\ldots,G_t).
$$
In fact, this follows from the standard mapping cone procedure that
relates the resolutions of $I, J$ and $\fc$.

Hence, keeping the notation above, we get an upper bound for the
Hilbert function of a relatively compressed algebra with respect to
$\fc$.

\begin{lemma} \label{lem-gen-upp-bound}
Let $A = R/I$ be a relatively compressed algebra with respect to $\fc$
and with socle degrees  $(s_1,\ldots,s_t)$. Let $G'_1,\ldots,G'_t \in R$
be general homogeneous forms of degree \linebreak
$e'-s_1,\ldots,e'-s_t$ and put $J' := \fc' + (G'_1,\ldots,G'_t)$.
Then we have for all integers $j$
$$
h_A(j) \leq h_{R/\fc'} (j) - h_{R/J'}(e'-j).
$$
\end{lemma}

\begin{proof}
By Liaison theory we have the formula
$$
h_A(j) = h_{R/\fc'} (j) - h_{R/J}(e'-j).
$$
Using Remark \ref{rem-cone} we get for all integers
$j$
$$
h_{R/J} (j) \geq h_{R/J'} (j)
$$
by the choice of the generators of $J'$. The claim follows.
\end{proof}

Note that the Hilbert function of $R/\fc'$ is determined by the
Hilbert function of $R/\fc$ and $s$. For example, if $c = n-1$ then
$$
h_{R/\fc'} (j) = h_{R/\fc}(j) - h_{R/\fc}(j-s-1)
$$
for all integers $j$.

If we start with $J'$  we get:

\begin{corollary} \label{cor-gen-min-b}
Let $F_1,\dots,F_{n-c} \in R$ be general forms of degree
$s+1$.
If the union of a minimal basis of $\fc$, $\{F_1,\dots,F_{n-c}\}$, and
$\{G'_1,\ldots,G'_t\}$ is
a minimal basis of $J'$ then we have equality in Lemma
\ref{lem-gen-upp-bound}, i.e.\
$$
h_A(j) = h_{R/\fc'} (j) - h_{R/J'}(e-j).
$$
\end{corollary}

\begin{proof}
The assumption on the minimal generators of $J'$ guarantees that $A =
R/I$, where $I := \fc' : J'$, has socle degrees $(s_1,\ldots,s_t)$. Hence
Lemma \ref{lem-gen-upp-bound}  shows that $A$ is compressed with respect to
$\fc$.
\end{proof}

\begin{remark} Note that by semicontinuity, the Hilbert function of
$J'$ does not
depend
    on the choice of the polynomials $G'_i$. It is determined by
    $R/\fc$ and the degrees $e'-s_1,\ldots,e'-s_t$ of
    the extra generators (since these were chosen generally).  In
general there are no
explicit formulas.  However, if $\fc$ is a general complete intersection then
Fr\"oberg's conjecture (see below) predicts the precise value of the Hilbert
function.  But independently of Fr\"oberg, if the assumption of Corollary
\ref{cor-gen-min-b}  is satisfied then  in  principle it allows us to
compute (at least on
the computer) the Hilbert  function of  a relatively compressed
algebra with respect
to
$\fc$  and with socle degrees  $(s_1,\ldots,s_t)$ from the given data.
\end{remark}

Thanks to a conjecture of Fr\"oberg (\cite{froberg}, page 120) all this
can be made (conjecturally) more explicit if $\fc$ is a general complete
intersection. In order to state this conjecture recall that the
Hilbert series of an algebra $A$ is the formal power series
$$
H_A (Z) := \sum_{j \geq 0} h_A (j) Z^j.
$$
We define for a power series ${\displaystyle \sum_{j \geq 0} a_j Z^j}$
with real coefficients
$$
\left | \sum_{j \geq 0} a_j Z^j \right | := \sum_{j \geq 0} b_j Z^j
$$
where
$$
b_j := \left \{ \begin{array}{ll}
a_j & \mif a_i\ge 0 \hbox{ for all $i\le j$} \\
0 & \mbox{otherwise}.
\end{array} \right.
$$

\begin{conjecture}[Fr\"oberg] Let $J \subset R$ be an ideal generated
    by general forms of degree $d_1,\ldots,d_r$.  Then the Hilbert series
    of $R/J$ is
$$
H_{R/J} (Z) =  \left | \frac{\prod_{i = 1}^r (1 - Z^{d_i})}{(1 - Z)^n}
\right |.
$$
\end{conjecture}

Note that it is easy to see that
$$
H_{R/J} (Z) \geq  \left | \frac{\prod_{i = 1}^r (1 - Z^{d_i})}{(1 - Z)^n}
\right |
$$
where the estimate compares the coefficients of same degree
powers.  Moreover, the conjecture was proved to be true for $n=2$ in
\cite{froberg}, for $n=3$ in \cite{anick}  and for arbitrary $n$ if $J$ is
a complete intersection or if $J$ is an almost complete
intersection \cite{stanley}.

Now we specialize to the case where $\fc \subset R$ is a complete
intersection of type
$(d_1,\ldots,d_c)$. Then the Castelnuovo-Mumford regularity of
$R/\fc$ is $e = d_1 +
\ldots + d_c - c$.

\begin{corollary}
Let $A = R/I$ be a relatively compressed algebra with respect to the
complete intersection $\fc$ and with socle degrees 
$(s_1,\ldots,s_t)$. Put
$$
e' = (n-c) s + d_1 + \ldots + d_c - c.
$$
Then
$$
H_A (Z) \leq \left | \frac{\prod_{i = 1}^c (1 - Z^{d_i}) \cdot (1 -
   Z^{s+1})^{n-c}}{(1 - Z)^n}
\right | - \left | \frac{\prod_{i = 1}^c (1 - Z^{d_i}) \cdot
(1 - Z^{s+1})^{n-c} \cdot
   \prod_{i=1}^t (1 - Z^{e'-s_i})}{(1 - Z)^n}
\right |.
$$
\end{corollary}

\begin{proof}
Using the notation of Lemma \ref{lem-gen-upp-bound} it is easy to see that
$$
H_{R/J'}(Z) \geq  \left | \frac{\prod_{i = 1}^c (1 - Z^{d_i}) \cdot
  (1 - Z^{s+1})^{n-c} \cdot
   \prod_{i=1}^t (1 - Z^{e'-s_i})}{(1 - Z)^n}
\right |.
$$
The claim follows.
\end{proof}

In fact, if $\fc$ is a {\em general}
complete intersection then we often expect equality.

\begin{conjecture} \label{conj hilb ser}
Let $J \subset R$ be the ideal generated by general forms
$G_1,\ldots,G_{c+t}$ of degree $d_1,\ldots,d_c,
e-s_1.\ldots,e-s_t$ and general forms $F_1,\ldots,F_{n-c} $ of degree
$s+1$. Assume that $\{G_1,\ldots,G_{c+t}, F_1,\ldots,F_{n-c}\}$ is a minimal
basis of $J$. Then the Hilbert series of a relatively compressed
algebra $A$  with respect to a {\rm general}
complete intersection of type $(d_1,\ldots,d_n)$  and with socle
degrees 
$(s_1,\ldots,s_t)$ is
$$
H_{A}(Z) = \left | \frac{\prod_{i = 1}^c (1 - Z^{d_i}) \cdot (1 -
   Z^{s+1})^{n-c}}{(1 - Z)^n}
\right | - \left | \frac{\prod_{i = 1}^c (1 - Z^{d_i}) \cdot
(1 - Z^{s+1})^{n-c} \cdot
   \prod_{i=1}^t (1 - Z^{e'-s_i})}{(1 - Z)^n}
\right |.
$$
\end{conjecture}

\begin{remark}
By Corollary \ref{cor-gen-min-b} this conjecture is true if and only
if Fr\"oberg's
conjecture is true in the corresponding case.
Note that granting Fr\"oberg's conjecture, the assumption of
Conjecture 
\ref{conj hilb ser} can be translated into a purely numerical condition
involving only the numbers
$d_1,\ldots,d_n,$ 
$e-s_1,\ldots,e-s_t$ in every specific example  though it seems difficult to
    make this explicit in general.
\end{remark}

\begin{remark} \label{char p remark}
Potentially, the Hilbert function of a relatively
compressed algebra with respect to a specific complete intersection
could differ from the one with respect to a {\em general} complete
intersection of the same type, and the result can change as the
characteristic varies. While it might not be the case over fields of
characteristic zero (it is an open question) that different complete
intersections yield different Hilbert functions for relatively compressed
algebras of the same type, this phenomenon does occur over fields of positive
characteristic.  The following example illustrates all these things.
\end{remark}

\begin{example}
We illustrate  the assertions of Remark \ref{char p remark}.  We will
consider the case of three variables, $R = k[x_1,x_2,x_3]$ (leaving open for
now the characteristic), and a complete intersection $\fc = (F_1,F_2,F_3)$
generated by three quartics.  This complete intersection has Hilbert function
\[
1 \ 3 \ 6 \ 10 \ 12 \ 12 \ 10 \ 6 \ 3 \ 1.
\]
We will be interested in relatively compressed Gorenstein algebras of socle
degree 8.  The expected Hilbert function for this algebra is
\[
1 \ 3 \ 6 \ 10 \ 12 \ 10 \ 6 \ 3 \ 1,
\]
and by \cite{Iarrobino-LNM} Theorem 4.16, this is achieved if $k$ has
characteristic zero and if, in addition, the complete intersection is either
the monomial complete intersection $(x_1^4, x_2^4, x_3^4)$, or general.

We now consider a field $k$ of characteristic 2.  A standard mapping cone
argument (as used also elsewhere in this paper) gives that the relatively
compressed Gorenstein algebra is linked via $\fc$ to an
ideal $I$ with four generators, and the fact that the desired socle degree of
the Gorenstein algebra is 8 yields that $I = (L,F_1,F_2,F_3)$, where $L$ is
a linear form.  The fact that the Gorenstein algebra has maximal Hilbert
function tells us that the linear form is general, and  that the Hilbert
function of $I$ is as small as possible among ideals with these generator
degrees.  Note also that in order to obtain the expected Hilbert function for
the relatively compressed Gorenstein algebra, the Hilbert function of $I$
must be
\[
1 \ 2 \ 3 \ 4 \ 2.
\]
We first consider the complete intersection $\fc = (x_1^4,x_2^4,x_3^4)$.
If $L = ax_1+bx_2+cx_3$ is a general linear form, then (because of the
characteristic) in fact
$I$ is a complete intersection, $I = (L, x_1^4,x_2^4)$.  Its Hilbert
function is
\[
1 \ 2 \ 3 \ 4 \ 3 \ 2 \ 1.
\]
This is not the required Hilbert function, so $R/\fc$ does not have a
relatively compressed Gorenstein algebra with the predicted Hilbert
function.  In fact, the Hilbert function of a relatively compressed
Gorenstein algebra with respect to this complete intersection is
\[
1 \ 3 \ 6 \ 9 \ 10 \ 9 \ 6 \ 3 \ 1.
\]
In fact, by studying the mapping cone for the link of $I$ via this complete
intersection, one sees that the relatively compressed Gorenstein algebra is
in fact itself a complete intersection of type $(3,4,4)$.

Now we consider what happens (still in characteristic 2) if we change the
generators.  We have verified using {\tt macaulay} that in characteristic 2,
if we change the complete intersection to
$(F_1,F_2,F_3)$ where
\[
\begin{array}{rcl}
F_1 & = &
x_1^4+x_1x_2^3+x_1^2x_2x_3+x_1^2x_3^2+x_1x_2x_3^2+x_1x_3^3+x_2x_3^3 \\
F_2 & = &
x_1^3x_2+x_1^2x_2x_3+x_1x_2^2x_3+x_2^3x_3+x_2^2x_3^2+x_2x_3^3+x_3^4 \\
F_3 & = &
x_1x_2^3+x_1^3x_3+x_1^2x_2x_3+x_1x_2^2x_3+x_2^3x_3+x_1^2x_3^2+x_2^2x_3^2+
x_1x_3^3
\end{array}
\]
and take $L =  x_2$, then the Hilbert function we get is the same as that
obtained in characteristic zero, the expected one, completing our
verification of the assertions in Remark
\ref{char p remark}.

Notice that similar behavior occurs when we look for relatively compressed
Gorenstein algebras with respect to any monomial complete intersection where
all the generators have the same degree and this degree is a multiple of the
characteristic.

The ideas in this example are directly related to the question of whether
$R/ \fc$ has the Weak Lefschetz Property -- cf.\ \cite{HMNW} Remark 2.9 and
Corollary 2.4.

\end{example}

\medskip

We now consider relatively compressed {\em level} algebras.  First recall
(cf.\
\cite{Iarrobino-compr}, \cite{FL}) that a {\em compressed} (in the
classical sense)
level Artinian algebra
$A$  with socle degree
$s$ and socle dimension 
$c$ has Hilbert function
\[
h_A (t) = \min \{ \dim R_t, c \cdot \dim R_{s-t} \}.
\]
(The idea is that from the left and from the right the function grows
as fast as
it theoretically can; so it is not hard to see that this is an upper
bound, but the
hard part is to show that this bound is actually achieved, and even
more there is an
irreducible parameter space for which the bound is achieved on a Zariski-open
subset.)  For example, when $n = 3, s = 10$ and $c = 3$ we get the
Hilbert function
of the compressed level algebra to be 1 3 6 10 15 21 28 30 18 9 3.

We now want to consider relatively compressed level Artinian algebras, and in
particular algebras that are relatively compressed in a complete
intersection.  So,
letting $\fc \subset R$ be a complete intersection, we consider relatively
compressed level quotients of $R/\fc$.  Occasionally such algebras
are compressed in
the classical sense, but not usually.  We are interested in both the Hilbert
function and the minimal free resolution of such algebras.

We first determine another upper bound for the Hilbert function of a
graded level Artinian algebra $A$ that is relatively compressed with
respect to a complete intersection $\fc \subset R$.

\begin{lemma} \label{first bound}
Let $A =
R/I$ be a graded level Artinian algebra of socle dimension 
$c$, socle  degree $s$ and
relatively compressed with respect to a complete intersection $\fc
\subset R$.  Then
\[
h_A(t) \leq \min \{ \dim (R/\fc)_t, c \cdot \dim (R/\fc)_{s-t} \}.
\]
\end{lemma}

\begin{proof}
We use the theory of inverse systems (cf.\ \cite{tonynotes},
\cite{Iarrobino-LNM}) and refer to those sources for
the necessary background.  The following short summary is taken from
\cite{GHMS} Chapter 5.  Let $S = k[y_1,\dots,y_n]$.  We consider the
action of
$R$ on
$S$ by differentiation: if $F \in S_j$ then $x_i \circ F =
(\frac{\partial }{\partial y_i}) F$.  There is an order-reversing
function from the ideals of $R$ to the $R$-submodules of $S$ defined by
\[
\phi_1 : \{ \hbox{ideals of $R$} \} \rightarrow \{ \hbox{$R$-submodules of
$S$} \}
\]
defined by
\[
\phi_1(I) = \{ F \in S \ | \ G \circ F = 0 \hbox{ for all } G \in I \}.
\]
This is a 1-1 correspondence, whose inverse $\phi_2$ is given by $\phi_2(M)
= \Ann_R(M)$.  We denote $\phi_1(I)$ by $I^{-1}$, called the {\em inverse
system} to $I$.  The pairing
\[
R_j \times S_j \rightarrow S_0 \cong k
\]
is perfect.  For a subspace $V$ of $R_j$ we write $V^\perp \subset S_j$ for
the annihilator of $V$ in this pairing.
If $F$ is an element of $[\fc_j]^{\perp}$ then the ideal
$J = \Ann(F)$ is a Gorenstein ideal containing $\fc$, of socle degree $j$.

Now we proceed by induction on $c$.  For $c=1$, we take $F \in
[\fc_s]^{\perp}$ and we
consider the Gorenstein graded algebra $R/\Ann(F)$.
Because $R/\Ann(F)$ is a  quotient of $R/\fc$, which is Gorenstein, we
clearly have
\begin{equation} \label{hilb fn bound}
h_A(t) \leq \min \{ \dim (R/\fc)_t, \dim (R/\fc)_{s-t} \}.
\end{equation}
Note that by \cite{Iarrobino-LNM}, Theorem 4.16, if $\fc$ and $F$ are
both  general (or if $\fc$ is a monomial complete intersection and $F$
is general) and if $k$ has characteristic zero (but see also
Remark \ref{characteristic remark}) then we have equality in (\ref{hilb fn
bound}) and
$R/\Ann(F)$ is relatively compressed.
To prove the general case we can choose independent elements $F_1,\dots,F_c
\in [\fc_s]^{\perp}$. Summing up the Hilbert functions of the Gorenstein
quotients
$R/\Ann(F_1),\dots,R/\Ann(F_c)$ of $R/\fc$, we get
\[
\dim (R/\Ann (F_1,\dots,F_c))_t \leq \sum_{i=1}^c \dim (R/\Ann (F_i))_t
\]
from which the result follows.
\end{proof}

It was noted in the proof above that when $c=1$ and $\fc$ and $F$ are 
general, and if the field has characteristic zero or satisfies a certain
numerical condition (see Remark \ref{characteristic remark}) then the
Gorenstein quotient that we get is relatively compressed with respect to
$\fc$, i.e.\ the inequality  (\ref{hilb fn bound}) is an equality.  It is
natural to ask if the same holds for $c \geq 2$.  We now present some
examples that show that the naive guess that the inequality in the statement
of Lemma \ref{first bound} is an equality even for ``generic'' choices, is
not always correct.

\begin{example} \label{orig 2.4}
Consider  the general complete intersection $\fc$
of type $(3,3,3)$ in 3 variables.  Its Hilbert function is 1 3 6 7 6
3 1.  Suppose
we want a level algebra $A = R/I$ with $s = 5$ and $c = 2$ that is relatively
compressed with respect to $R/\fc$.
    One would ``expect" that its Hilbert function  would be 1 3 6 7 6 2.

First, we note that there is an algebra with this Hilbert function, which is a
quotient
of that complete intersection, but it {\em must} have a ghost term making it
not be level.  The reason comes from liaison theory.  Let $J$ be the
residual to $I$
with respect to the complete intersection $\fc$.  Then $R/J$ would have Hilbert
function 1 1 (cf.\ \cite{migbook} Corollary 5.2.19),  so its resolution begins
\[
\cdots \rightarrow R(-2)
\oplus R(-3)^2
\to R(-1)^2 \oplus R(-2) \to J \to 0.
\]
   Since the complete intersection $\fc$ does not contain a quadric,
the mapping cone
procedure (cf.\ \cite{migbook} Proposition 5.2.10) does not split off
any summands
corresponding to generators of $J$, so the ideal $I$ has resolution that ends
\[
0 \rightarrow R(-8)^2  \oplus R(-7)  \rightarrow R(-7) \oplus R(-6)^5
\rightarrow \dots
\]
and so is not level.

Second, we show that there is a relatively compressed algebra that has Hilbert
function 1 3 6 7 5 2.  Indeed, it can be obtained as the residual of a complete
intersection of type $(1,1,3)$ inside  complete  intersection of type
$(3,3,3)$.

Finally, we observe that all this can even be done at the level of
points in $\PP^3$
by  lifting the Artinian ideals to ideals of sets of points: simply
start with a set
of three points on a line in $\PP^3$ and link using three cubics.
\end{example}

\begin{example} Now we consider ideals in the polynomial ring $R$ with 4
variables.  Let $\fc$ be a general complete intersection of type
(3,3,3,3).  Its Hilbert
function is
\[
1 \ 4 \ 10 \ 16 \ 19 \ 16 \ 10 \ 4 \ 1.
\]
    Suppose we look for a relatively compressed
level algebra with s = 7 and c = 2.  A first guess could be that the correct
Hilbert function should be 1 4 10 16 19 16 8 2.  If this were true,
the residual $J$
in the complete intersection $(3,3,3,3)$ would have Hilbert function
\linebreak 1 2
2.  But such
$J$  has (at least) generators of degrees 1, 1, 2, 3, 3, so the mapping cone
procedure shows that
$R/I$ has Cohen-Macaulay type (at least) 3, hence it is not a level
algebra. Again
there is a ghost term.  Note that generically we get Cohen-Macaulay
type exactly 3.

On the other hand, we  can again construct a relatively compressed
level algebra with $s = 7$, $c = 2$, and  Hilbert
function 1 4 10 16 19 16 7 2.  This is done by starting with an algebra
with Hilbert function 1 2 3.  Its generators will be of degree 1, 1, 3, 3, 3, 3
so the complete intersection will split off all the terms
corresponding to the cubic
generators, leaving a residual that is a level algebra.
\end{example}

\begin{example} \label{orig 2.6}
   This time we will even compute the graded Betti
numbers, not just  the ``surprising'' Hilbert function. We work over
a polynomial
ring with 3 variables.  Recall that in this case we know Fr\"oberg's
conjecture to
hold, thanks to work of Anick \cite{anick},  which gives the Hilbert
function of an
ideal generated by general forms of any prescribed degrees (as
illustrated below).

Consider a level algebra with $s = 7$, $c = 2$ that is relatively
compressed with
respect to the general complete intersection $\fc$  of type $(4, 4, 4)$. It is
constructed as follows.   Adjoining two general forms of degree 2 to
the general
complete intersection $\fc$ of type $(4, 4, 4)$ we get (using Anick's result)
successively the Hilbert functions
$$
\begin{array}{ccccccccccccc}
&&1& 3& 6& 10& 12& 12& 10& 6& 3 & 1\\
&&1& 3& 5&  7&  6&  2\\
&& 1& 3& 4&  4&  1.
\end{array}
$$
The residual with respect to $\fc$ provides the desired algebra $A$. It has
Hilbert function
$$
1 \ 3 \ 6 \ 10 \ 12 \ 11 \  6 \ 2,
$$
where one might have expected in degree 5 a 12 rather than 11. The
reason that there
is no level algebra with a 12 rather than an 11 in degree 5, that is
a quotient of
$R/\fc$, is precisely that such an algebra would be residual to an algebra with
Hilbert function \linebreak 1, 3, 4, 4, which has too many generators
to allow the 
residual to be level! Note however that there is a level algebra $R/I$ with
Hilbert function $1, 3, 6, 10, 12, 12,  6, 2$ if we do not require that
$I$ contains a regular sequence of type $(4, 4, 4)$. It can be
constructed as Artinian quotient of the coordinate ring of 12 points
in $\PP^2$ using \cite{GHMS}, Proposition 7.1. 

In order to compute the graded Betti numbers of $A$ we start with the
set $X$ of 3
general points in $\PP^3$. Their resolution has the shape
$$
0 \to  R^2(-4) \to R^5(-3) \to R(-1) \oplus R^3(-2) \to I(X) \to 0.
$$

Linking by a complete intersection of type $(2,2,4)$ we get a
residual $J$ whose
Betti numbers read as
$$
0 \to R(-7) \oplus  R(-6) \to R^5(-5) \oplus R(-4) \to R^3(-4) \oplus
R^2(-2) \to J
\to 0
$$
because the two generators of degree 2 split off while the Koszul
ghost term does
not split off.

Since $J$ contains a complete intersection of type (2,2,4), it certainly
contains one of type (4,4,4) as well.  Linking again, this time by a complete
intersection of type
$(4,4,4)$, we get an algebra $A$ as above. After splitting off the
three quartics
the mapping cone procedure provides that its  minimal free resolution
has the form
$$
0 \to R^2(-10) \to R(-8) \oplus R^5(-7) \to R(-6) \oplus R(-5) \oplus
R^3(-4) \to J
\to 0.
$$
\end{example}


\section{Minimal free resolution of relatively compressed
Gorenstein Artinian
   algebras} \label{sec-Gor}

In the previous section we discussed what the ``expected'' behavior
should be for
the Hilbert function of a relatively compressed level algebra, and how this
is sometimes not achieved.  We now begin our study of the following
problem:

   \begin{problem} To determine the ``generic" graded Betti
numbers (or, equivalently, the minimal free $R$-resolution) of
Artinian level graded algebras of embedding dimension $n$, socle
degree $s$ and socle dimension $c$ and relatively compressed with respect to a
general complete intersection $\fa \subset R $.
\end{problem}

In this section we consider the minimal free $R$-resolution of
Gorenstein Artinian
graded algebras of embedding dimension $n$ that are relatively compressed
with respect to the ideal $\fa \subset k[x_1,...,x_n]$ of a general complete
intersection of type $(d_1,...,d_r)$, $r\leq n$.   We consider the case of
even socle degree and odd socle degree separately.  For even socle degree we
completely determine the minimal free $R$-resolution.  We also show that all
redundant (``ghost'') terms that appear are due to Koszul syzygies
(or are forced by
duality from Koszul syzygies).
For odd socle degree we do not have quite as clean a statement, but we show
in Example \ref{ghosts} that this is the best we could hope for.

We begin with even socle degree.  We will see that once we fix $n$, $t$ and
$(d_1,...,d_r)$ {\em all} have the same graded Betti numbers (there
is no need for a
``generic'' choice). This was known only for
compressed level algebras of even socle degree:

\begin{proposition}\label{res-comp-gor}  Let $A$ be a  compressed
   Gorenstein Artinian graded algebra  of embedding dimension
$n$ and socle degree $2t$.   Then $A$ has a minimal free
   $R$-resolution of the following type:
\[
0 \rightarrow R(-2t-n)\rightarrow R(-t-n+1)^{\alpha_{n-1}}
\rightarrow ... \rightarrow R(-t-p)^{\alpha _{p}} \rightarrow
\]

   \[.... \rightarrow
    R(-t-2)^{\alpha_2} \rightarrow
   R(-t-1)^{\alpha_1} \rightarrow  R \rightarrow A \rightarrow 0
\]
where
\[
\alpha_{i}={t+i-1 \choose i-1}{t+n \choose n-i}-{t-1+n-i \choose
n-i}{t-1+n\choose i-1}
\]
   for $i=1,...,n-1$.
\end{proposition}

\begin{proof}
Because the socle degree is even, in fact $A$ is a so-called {\em extremely
compressed} Artin level algebra.  Then the result follows from \cite{boij2}
Proposition 3.6 or \cite{Iarrobino-compr} Proposition 4.1.
\end{proof}

\begin{remark}
The above result is also a special case of \cite{MN3}, Theorem 8.14.
This latter result
has the extra hypothesis that $A$ has the Weak Lefschetz Property (i.e.\
that multiplication by a general linear form, from any component to
the next, has
maximal rank).  However, it was noted in
\cite{MMR2} Remark 3.6(c) that in the situation of compressed
Gorenstein algebras of
even socle degree, this property is automatically satisfied.

Note also that the formula for $\alpha_i$ given above is not
presented in the same
way as it is in \cite{MN3} Theorem 8.14, but a calculation shows that they are
equivalent.
\end{remark}

\vskip 2mm
   From now on, when we say that $A$ is a Gorenstein Artinian graded
algebra of embedding dimension $n$, even socle degree $2t$ and
relatively compressed with respect to a general complete intersection
ideal $\fa \subset k[x_1,...,x_n]$ of   type $d_1\le \cdots \le
d_r$ we will assume without loss of generality that $d_r \le t$;
otherwise $A$ is a
Gorenstein Artinian graded algebra of embedding dimension $n$, even
socle degree $2t$ and relatively compressed with respect to a
general complete intersection  ideal $\fb \subset k[x_1,...,x_n]$ of type
$d_1\le \cdots \le d_j$ where $d_j=\max_{1\le i \le r} \{ d_i\le t
\}.$

\vskip 4mm We first fix some notation that we will use from now on.

\begin{notation}

\item[(1)] Given a complete intersection ideal  $\fa  =
(G_1,\dots,G_r)\subset R=k[x_1,...,x_n]$ with $r \le n$ and
$d_1=deg(G_1) \leq \dots \leq d_r=deg(G_r)$, we denote by
$K_i(d_1,\cdots ,d_r)$ (or, simply, $K_i(\underline{d})$ if
$\underline{d}=(d_1,\cdots ,d_r)$) the $i$-th module of syzygies
of $R/\fa $. So, we have
\[
K_i(d_1,\cdots ,d_r)= K_i(\underline{d}) := \bigwedge^i \left (
\bigoplus_{i=1}^r R(-d_i) \right )
\]
and the minimal free $R$-resolution of $R/\fa $:
\[
0 \rightarrow K_r (\underline{d})\rightarrow
K_{r-1}(\underline{d}) \rightarrow \cdots \rightarrow
K_2(\underline{d}) \rightarrow K_1 (\underline{d})\rightarrow R
\rightarrow R/\fa \rightarrow 0.
\]

\item[(2)] For any free $R$-module $F=\oplus _{t\in
\ZZ}R(-t)^{\alpha _t}$ and any integer $y \in \ZZ$, we set
\[ F^{\le y}:=\oplus _{t\le y
}R(-t)^{\alpha _t}. 
\]

\end{notation}

\begin{theorem}\label{res-rel-comp-gorcodim-n}
   Let $A=R/I$ be a
    Gorenstein Artinian graded algebra  of embedding dimension
$n$ and socle degree $2t$, where $R = k[x_1,\dots,x_n]$ and $k$ has
characteristic zero.  Assume that
$A$ is relatively compressed with respect to a general complete intersection
ideal
$\fa =(G_1, \cdots ,G_r)$, $r\le n$ and $\deg(G_1)=d_1$,\dots,
$\deg(G_r)=d_r$. Set $\underline{d}=(d_1, \cdots , d_r)$. Then,
$A$ has a minimal free
    $R$-resolution of the following type:

\[
0 \rightarrow R(-2t-n)  \rightarrow R(-t-n+1)^{\alpha
_{n-1}(\underline{d},n,t)} \oplus K_{n-1}(\underline{d})^{\le
t+n-2} \oplus K_{1}(\underline{d})^{\vee }(-2t-n) \rightarrow
\]
\[R(-t-n+2)^{\alpha
_{n-2}(\underline{d},n,t)} \oplus K_{n-2}(\underline{d})^{\le
t+n-3} \oplus (K_{2}(\underline{d})^{\le t+1})^{\vee }(-2t-n)
\rightarrow \cdots \rightarrow
\]
\[\rightarrow R(-t-2)^{\alpha _{2}(\underline{d},n,t)} \oplus
K_{2}(\underline{d})^{\le t+1} \oplus (K_{n-2}(\underline{d})^{\le
t+n-3})^{\vee }(-2t-n) \rightarrow
\]
\[\oplus
_{j=1}^{r}R(-d_j) \oplus R(-t-1)^{\alpha
_{1}(\underline{d},n,t)}\oplus (K_{n-1}(\underline{d})^{\le
t+n-2})^{\vee }(-2t-n) \rightarrow R \rightarrow A \rightarrow 0
\]
where \[\alpha _i(\underline{d},n,t)= \alpha
_{n-i}(\underline{d},n,t)\mbox{ for } i=1, \cdots n-1,
\]
and $\alpha _i(\underline{d},n,t)$ is completely determined by the
Hilbert function of $A$.

\end{theorem}
\begin{proof}
Since $\fa \subset I$ and $A$ is relatively compressed with respect to
$\fa$, we have $ h_A(\nu) = \min \{ \dim
(R/\fa)_\nu, \dim (R/\fa)_{2t-\nu} \}$ (thanks to \cite{Iarrobino-LNM},
Theorem 4.16).  Thus $\fa _\nu =I_\nu$ for all $\nu \le t$.  We
deduce  that
$A$ has a minimal free
$R$-resolution of the following type

\vskip 2mm
\begin{equation} \label{initial}
0 \rightarrow R(-2t-n) \rightarrow F_{n-1} \rightarrow F_{n-2}
\rightarrow \cdots \rightarrow F_{2} \rightarrow F_{1}\rightarrow
R \rightarrow R/I \rightarrow 0
\end{equation}

\vskip 2mm \noindent where \[ F_i= \oplus _{m\ge t+i}R(-m)^{\alpha
_{i,m}(\underline{d},n,t)}\oplus K_{i}(\underline{d})^{\le
t+i-1},\]

\vskip 2mm \noindent$\underline{d}=(d_1,\cdots ,d_r)$,
$K_i(\underline{d})=K_i(d_1,\cdots ,d_r)$ and $\alpha
_{1,t+1}(\underline{d},n,t)$ is completely determined by the
Hilbert function of $A$. More precisely, $\alpha
_{1,t+1}(\underline{d},n,t)=\dim (R/ \fa )_{t+1}-\dim (A)_{t+1}$.

\vskip 2mm Recall that the minimal free $R$-resolution of a
Gorenstein Artinian algebra is self-dual  (up to shift). Dualizing
(\ref{initial}) and twisting by $R(-2t-n)$, we get
\[
\begin{array}{c}
0 \rightarrow R(-2t-n) \rightarrow  F_1^{\vee}(-2t-n)\rightarrow
F_2^{\vee}(-2t-n)\rightarrow \cdots \hbox{\hskip 1.4in}
\\ \vbox{ \vskip .2in}
\hfill F_{n-2}^{\vee}(-2t-n)\rightarrow F_{n-1}^{\vee}(-2t-n)\rightarrow
R \rightarrow A=R/I \rightarrow 0.
\end{array}
\]

Therefore, for all $1\le i \le n-1$, we have the isomorphism
\begin{equation} \label{isom}
\begin{array}{rcl}
  F_i & = & \oplus _{m\ge t+i}R(-m)^{\alpha_{i,m}(\underline{d},n,t)}\oplus
K_{i}(\underline{d})^{\le t+i-1} \\ \vbox{ \vskip .2in}
& \cong & F_{n-i}^{\vee}(-2t-n) \\ \vbox{ \vskip .2in}
& = & \oplus _{m\ge t+n-i}R(m-2t-n)^{\alpha
_{n-i,m}(\underline{d},n,t)}\oplus (K_{n-i}(\underline{d})^{\le
t+n-i-1})^{\vee}(-2t-n).
\end{array}
\end{equation}
We first consider $i=1$.  We rewrite the third line of (\ref{isom}) as
follows:
\[
\begin{array}{c}
\left [
R(-t-1)^{\alpha_{n-1,t+n-1}(\underline{d},n,t)} \oplus
R(-t)^{\alpha_{n-1,t+n} (\underline{d},n,t)} \oplus \cdots \right ] \oplus \\
\vbox{ \vskip .2in}
\left [ (K_{n-1}(\underline{d})^{\leq
t+n-2})^\vee(-2t-n) \right ]
\end{array}
\]
and observe that each summand on the first line is of the form $R(-i)$ for
some $i \leq t+1$, while each summand on the  second line is of the form
$R(-i)$ for some $i \geq t+2$.
It follows that
$$\oplus _{m\ge t+n-1}R(m-2t-n)^{\alpha
_{n-1,m}(\underline{d},n,t)} =K_{1}(\underline{d})^{\le t}\oplus
R(-t-1)^{\alpha _{1,t+1}(\underline{d},n,t)}, \mbox{ and }$$
$$(K_{n-1}(\underline{d})^{\le
t+n-2})^{\vee}(-2t-n)=\oplus _{m \ge t+2}R(-m)^{\alpha
_{1,m}(\underline{d},n,t)}$$ and we conclude that

$$F_1=K_{1}(\underline{d})^{\le t}\oplus R(-t-1)^{\alpha
_{1}(\underline{d},n,t)}\oplus (K_{n-1}(\underline{d})^{\le
t+n-2})^{\vee}(-2t-n); \mbox{ and }$$

$$F_{n-1}=(K_{1}(\underline{d})^{\le t})^{\vee}(-2t-n)\oplus R(-t-n+1)^{\alpha
_{1}(\underline{d},n,t)}\oplus K_{n-1}(\underline{d})^{\le
t+n-2}$$

\noindent where  $\alpha _{1}(\underline{d},n,t):=\alpha
_{1,t+1}(\underline{d},n,t)$.

\vskip 2mm Substituting $F_1$ and $F_{n-1}$ in the exact sequence
(\ref{initial}) and using again  the Hilbert function of $A$, we
determine $\alpha _{2,t+2}(\underline{d},n,t)$. Moreover, an
analogous numerical analysis taking into account that $F_2\cong
F_{n-2}^{\vee }(-2t-n)$ gives us

$$F_2=K_{2}(\underline{d})^{\le t+1}\oplus R(-t-2)^{\alpha
_{2}(\underline{d},n,t)}\oplus (K_{n-2}(\underline{d})^{\le
t+n-3})^{\vee}(-2t-n); \mbox{ and }$$

$$F_{n-2}=(K_{2}(\underline{d})^{\le t})^{\vee}(-2t-n)\oplus R(-t-n+2)^{\alpha
_{2}(\underline{d},n,t)}\oplus K_{n-2}(\underline{d})^{\le
t+n-3}$$

\noindent where  $\alpha _{2}(\underline{d},n,t):=\alpha
_{2,t+2}(\underline{d},n,t)$.

\vskip 2mm Going on and using the isomorphism $F_i\cong
F_{n-i}^{\vee }(-2t-n)$ for all $i=1, \cdots , n-1$, we obtain
that

$$F_i=K_{i}(\underline{d})^{\le t+i-1}\oplus R(-t-i)^{\alpha
_{i}(\underline{d},n,t)}\oplus (K_{n-i}(\underline{d})^{\le
t+n-i-1})^{\vee}(-2t-n)$$

\noindent where  $\alpha _{i}(\underline{d},n,t)=\alpha
_{n-i}(\underline{d},n,t)$ for all $i=1, \cdots, n-1$ and $\alpha
_{i}(\underline{d},n,t)$ is determined by the Hilbert function of
$A$.
\end{proof}

\begin{remark}\label{characteristic remark}
The equality in the first line of the proof of Theorem
\ref{res-rel-comp-gorcodim-n} follows from the assumption that the complete
intersection is general, as well as the assumption on the characteristic.
We remark that the assumptions can be weakened somewhat.

As we saw in Section \ref{sec-exist}, the fact that the Gorenstein algebra
$A$ has the expected Hilbert function (i.e.\ that it is relatively
compressed with respect to $\fa$) is directly related to the Fr\"oberg
conjecture, in this case for $n+1$ forms.  This in turn is equivalent to the
so-called {\em Maximal Rank Property}, namely that $R/\fa$ have the
property that for any $d$ and any $i$, a general form $F$ of degree $d$
induces a map of maximal rank from $(R/\fa)_i$ to $(R/\fa)_{i+d}$.  And this
follows from the Strong Lefschetz Property.  Now, it was shown in
\cite{Iarrobino-LNM} (based on the proof of Fr\"oberg's conjecture for
$n+1$ forms in $n$ variables in \cite{stanley} and \cite{watanabe})
that all of these hold for a monomial complete 
intersection (hence for a general complete intersection) provided that
either $k$ has characteristic zero or that $k$ has characteristic $p$,
assuming that certain numerical conditions hold.  More precisely, they
assume that either char$(k) = 0$ or else that char$(k) > j$,  $\fa =
(f_1,\dots,f_a), a \leq n$, $f_i = x_i^{d_i}$, and that there exist
nonnegative integers $t, u, v$ with $j = u+v$, and if $a=n$ then $j \leq
\sum d_i -n$. Hence Theorem \ref{res-rel-comp-gorcodim-n} is  also true
with these assumptions on the characteristic.

Note that in fact it is unknown if the Strong Lefschetz Property  holds for
{\em all} complete intersections, even in characteristic
zero. However, it is true that over fields of characteristic zero all
complete intersections in 3 variables have the Weak Lefschetz
Property,  due to \cite{HMNW}. 
\end{remark}

\begin{example} \label{examp1}Let $A$ be a general Gorenstein Artinian
graded  algebra of embedding
dimension 4, socle degree 10 and relatively compressed with
respect to the ideal $\fa =(F,G,H)$ of a complete intersection set
of points $P\subset \PP^3$ of type (3,3,4). The $h$-vector of $A$
is
\[
1 \ 4 \ 10 \ 18 \ 26 \ 32 \ 26 \ 18 \ 10 \  4 \ 1
\]
and the ``expected" minimal free $R$-resolution is
\[
\begin{array}{c}
0 \rightarrow R(-14) \rightarrow R(-8)^{9}\oplus R(-11)^2\oplus
R(-10) \rightarrow R(-7)^{20}\oplus R(-8)\oplus R(-6) \hbox{\hskip
1in}
\\ \\
\hbox{\hskip 2in} \rightarrow R(-6)^{9}\oplus R(-3)^2\oplus R(-4)
\rightarrow R \rightarrow A \rightarrow 0.
\end{array}
\]
It is known that such algebras exist \cite{Iarrobino-LNM}, and the precise
resolution comes from
Theorem \ref{res-rel-comp-gorcodim-n}.  However, to illustrate our
technique from
the previous section we explicitly  construct it. To this end, we consider a
subset $X\subset P \subset \PP^3$ of  32  points  of $P$ that truncate the
Hilbert function. The
$h$-vector of $X$ is thus
\[
1 \ 3 \ 6 \ 8 \ 8 \ 6 \ 0 \
\]
    and $I(X)$ has a minimal free $R$-resolution of the
following type
\[
0 \rightarrow R(-8)^6 \rightarrow R(-7)^{10}\oplus R(-6)
\rightarrow R(-6)^{3}\oplus R(-4)\oplus R(-3)^2 \rightarrow R
\rightarrow R/I(X) \rightarrow 0
\]
(as can  be verified, for example, by linkage).

The canonical module $\omega_X$ of $R/I(X)$ can be embedded as an
ideal $\omega_X (-10)\subset R/I(X)$ of initial degree 6 and we
have a short exact sequence
\[
0 \rightarrow \omega_X (-10)\rightarrow R/I(X) \rightarrow A
\rightarrow 0
\]
where $A$ is a Gorenstein Artinian graded algebra of codimension
4, socle degree 10, $h$-vector
\[
1 \ 4 \ 10 \ 18 \ 26 \ 32 \ 26 \ 18 \ 10 \  4 \ 1.
\]
So, it is relatively compressed with respect to $\fa =(F,G,H)$,
$\deg(F)=\deg(G)=3$ and $\deg(H)=4$. Moreover, applying once more the
mapping cone process we get that $A$ has the following minimal
free $R$-resolution

\[
\begin{array}{c}
0 \rightarrow R(-14) \rightarrow R(-8)^{9}\oplus R(-11)^2\oplus
R(-10) \rightarrow R(-7)^{20}\oplus R(-8)\oplus R(-6) \hbox{\hskip
1in}
\\ \\
\hbox{\hskip 2in} \rightarrow R(-6)^{9}\oplus R(-3)^2\oplus R(-4)
\rightarrow R \rightarrow A \rightarrow 0.
\end{array}
\]
So, it has the expected minimal free $R$-resolution in the sense
that the graded Betti numbers are the smallest consistent with the
Hilbert function of such relatively compressed Gorenstein algebra.

We see that the summand $R(-6)$ does not split off because it is a
Koszul relation
among the two cubic generators.  The summand $R(-8)$ does {\em not}
correspond to a
Koszul syzygy.  However, it is forced by the self duality property
(up to twist) of
the minimal free $R$-resolution of an Artinian Gorenstein graded
algebra. So, in
the minimal free resolution, the ghost terms are forced to be there
(directly and
then indirectly) by the Koszul relations among the generators.
\end{example}

\begin{remark}
Theorem \ref{res-rel-comp-gorcodim-n} shows that the observation at
the end of the
last example holds in general: in the minimal free resolution of an Artinian
Gorenstein algebra of even socle degree, relatively compressed with
respect to a
general complete intersection, the only ghost terms that appear are
forced to be
there by Koszul relations among the generators, or by duality because
of such Koszul
relations.
\end{remark}

   A more difficult situation is when the Artinian
Gorenstein graded algebra has odd socle degree.  The technique of the previous
section only gives a partial result, not the precise resolution.
Indeed, it is no
longer true that the Hilbert function alone determines the graded
Betti numbers.
(For instance see \cite{MMR2} Example 3.12.)

\begin{theorem}\label{res-rel-comp-gor-oddsocle}
   Let $A=R/I$ be a
    Gorenstein Artinian graded algebra  of embedding dimension
$n$ and socle degree $2t+1$, where $R = k[x_1,\dots,x_n]$ and $k$ has
characteristic zero.  Assume that $A$ is relatively
compressed with respect to a general complete intersection ideal
$\fa =(G_1, \cdots ,G_r)$, $r\le n$ and $deg(G_1)=d_1, \cdots
, deg(G_r)=d_r$. Set $\underline{d}=(d_1, \cdots , d_r)$. Then,
$A$ has a minimal free
    $R$-resolution of the following type:

\[
0 \rightarrow R(-2t-1-n)\rightarrow F_{n-1} \rightarrow \cdots
\rightarrow F_2 \rightarrow F_1\rightarrow R \rightarrow A
\rightarrow 0 \]

\noindent where $F_i\cong F_{n-i}^{\vee }(-2t-1-n)$ for all $i=1,
\cdots , n-1$. Moreover, if $n$ is even (say, $n=2p$), then

\[
\begin{array}{c}
F_i=R(-t-i-1)^{y_{i+1}} \oplus R(-t-i)^{\alpha
_{i}(\underline{d},n,t)+y_i} \oplus  \hbox{\hskip 3in} \\
\hbox{\hskip 1in} K_{i}(\underline{d})^{\le t+i-1} \oplus
(K_{n-i}(\underline{d})^{\le t+n-i-1})^{\vee }(-2t-1-n) \mbox{ for
} i=2,\cdots , p-1;\end{array}
\]

$$F_1=\oplus
_{j=1}^{r}R(-d_j)\oplus (K_{n-1}(\underline{d})^{\le t+n-2})^{\vee
}(-2t-1-n) \oplus R(-t-1)^{\alpha _{1}(\underline{d},n,t)}\oplus
R(-t-2)^{y_2}; \mbox{ and }$$

$$F_p=R(-t-p-1)^{\alpha _{p}(\underline{d},n,t)+y_{p}} \oplus
R(-t-p)^{\alpha _{p}(\underline{d},n,t)+y_p} \oplus
K_{p}(\underline{d})^{\le t+p-1} \oplus (K_{p}(\underline{d})^{\le
t+p-1})^{\vee }(-2t-1-n)$$

\vskip 2mm \noindent where $\alpha _i(\underline{d},n,t)$, $i=1,
\cdots p$, is completely determined by the Hilbert function of
$A$.

\vskip 2mm If $n$ is odd (say, $n=2p+1$), then

$$F_1=\oplus
_{j=1}^{r}R(-d_j)\oplus (K_{n-1}(\underline{d})^{\le t+n-2})^{\vee
}(-2t-1-n) \oplus R(-t-1)^{\alpha _{1}(\underline{d},n,t)}\oplus
R(-t-2)^{y_2}; \mbox{ and }$$
\[
\begin{array}{c}
F_i=R(-t-i-1)^{y_{i+1}} \oplus R(-t-i)^{\alpha
_{i}(\underline{d},n,t)+y_i} \oplus  \hbox{\hskip 3in} \\
\hbox{\hskip 1in} K_{i}(\underline{d})^{\le t+i-1} \oplus
(K_{n-i}(\underline{d})^{\le t+n-i-1})^{\vee }(-2t-1-n) \mbox{ for
} i=2,\cdots , p\end{array}
\]

\vskip 2mm \noindent where $\alpha _i(\underline{d},n,t)$, $i=1,
\cdots p$, is completely determined by the Hilbert function of
$A$.

\end{theorem}

\begin{proof} Analogous to the proof of Theorem
\ref{res-rel-comp-gorcodim-n}.  See also Remark \ref{characteristic remark}
concerning the assumption on the characteristic.
\end{proof}

\begin{example} \label{ghosts}
This example shows that in fact a relatively compressed Gorenstein
algebra of odd
socle degree can have ghost terms that are not forced by Koszul
relations among the
generators, even taking duality into account.
   It was verified by {\tt macaulay}, but the calculations
can be done by hand as well. Let
$I$ be an ideal of general forms of degrees 4,4,4,4,11 in four
variables.  Link $I$
using a general complete intersection of type (4,4,4,11).  The residual is a
Gorenstein ideal $G$, which is relatively compressed since $I$ is an ideal of
general forms.  The Hilbert functions are

\bigskip

\begin{tabular}{l|ccccccccccccccccccccccccccccccccccccccccccc}
deg & 0 & 1 & 2 & 3 & 4 & 5 & 6 & 7 & 8 & 9 & 10 & 11 & 12 & 13 & 14
& 15 & 16 & 17
& 18 & 19 \\
CI & 1 & 4 & 10 & 20 & 32 & 44 & 54 & 60 & 63 & 64 & 64 & 63 & 60 & 54 &
44 & 32 & 20 & 10 & 4 & 1 \\\hline
G &  1 & 4 & 10 & 20 & 32 & 44 & 54 & 60 & 60 & 54 & 44 & 32 & 20 & 10 & 4 & 1
\end{tabular}

\bigskip

The Betti diagram of G is

\begin{verbatim}
total:      1     7    12     7     1
--------------------------------------
      0:      1     -     -     -     -
      1:      -     -     -     -     -
      2:      -     -     -     -     -
      3:      -     3     -     -     -
      4:      -     -     -     -     -
      5:      -     -     -     -     -
      6:      -     -     3     -     -
      7:      -     3     3     1     -
      8:      -     1     3     3     -
      9:      -     -     3     -     -
     10:      -     -     -     -     -
     11:      -     -     -     -     -
     12:      -     -     -     3     -
     13:      -     -     -     -     -
     14:      -     -     -     -     -
     15:      -     -     -     -     1
\end{verbatim}

The three copies of $R(-8)$ in the second free module represent
Koszul syzygies,
hence the corresponding ``ghost'' terms are forced by Koszul
relations.  However,
the copy of $R(-9)$ in the second free module does not come from
Koszul relations
among generators.  It  illustrates that things are very different in
the case of odd
socle degree.
\end{example}

\begin{remark}
The two preceding examples illustrate the difference between the case
of even socle
degree and odd socle degree.  In fact, Theorem
\ref{res-rel-comp-gorcodim-n} shows
that there are no ``ghost'' terms in the minimal free resolution of a
relatively
compressed Gorenstein Artinian  algebra of even socle degree, relatively
compressed with respect to a general complete intersection, apart from
those corresponding to Koszul relations or forced from the Koszul ones by
duality.
\end{remark}

Now we introduce a new technique, which gives the precise resolution
at least in a
special case.

\begin{example} Let $A$ be a general Gorenstein Artinian graded algebra of
embedding dimension 4, socle degree 9 and relatively compressed with respect
to the ideal $\fa =(F)$ of a smooth quadric $Q=V(F)\subset \PP^3$. The
$h$-vector of $A$ is
\[
1 \ 4 \ 9 \ 16 \ 25 \ 25 \ 16 \ 9 \ 4 \ 1
\]
and the expected minimal free $R$-resolution is
\[
0 \rightarrow R(-13) \rightarrow R(-11)\oplus R(-8)^{11}
\rightarrow R(-6)^{11}\oplus R(-7)^{11} \rightarrow R(-2)\oplus
R(-5)^{11} \rightarrow R \rightarrow A \rightarrow 0.
\]

Let us explicitly  construct such an algebra. To
this end, we consider 30 general points $X\subset \PP^3$ on the quadric
$Q\subset
\PP^3$. The $h$-vector of $X$ is 1 3 5 7 9 5 and $I(X)$ has a minimal free
$R$-resolution of the following type \cite{GMR}:
\[
0 \rightarrow R(-8)^5 \rightarrow R(-6)^{5}\oplus R(-7)^6
\rightarrow R(-2)\oplus R(-5)^{6} \rightarrow R \rightarrow R/I(X)
\rightarrow 0.
\]
By \cite{boij}; Theorem 3.2, the canonical module $\omega_X$ of
$R/I(X)$ can be embedded as an ideal $\omega_X (-9)\subset R/I(X)$
of initial degree $5$ and we have a short exact sequence
\[
0\rightarrow \omega_X (-9)\rightarrow R/I(X) \rightarrow A
\rightarrow 0
\]
where $A$ is a Gorenstein Artinian graded algebra of embedding
dimension 4, socle degree 9, $h$-vector
\[
1 \ 4 \ 9 \ 16 \ 25 \ 25 \ 16 \ 9 \ 4 \ 1
\]
and relatively compressed with respect to $\fa =(F)$. Moreover
applying the mapping cone process we get that $A$ has the
following minimal free $R$-resolution
\[
0 \rightarrow R(-13) \rightarrow R(-11)\oplus R(-8)^{11}
\rightarrow R(-6)^{11}\oplus R(-7)^{11} \rightarrow R(-2)\oplus
R(-5)^{11} \rightarrow R \rightarrow A \rightarrow 0.
\]

\end{example}

\vskip 4mm The following result from \cite{GMR} is crucial if we
want to generalize the above example.

\begin{proposition}\label{MRC} Let $X$ be a set of $N$ general points on a
smooth quadric $Q\subset \PP^3$. Write $N=i^2+h$ with $0<h\le
2i+1$. Then, the $h$-vector of $X$ is
\[ 1 \ 3 \ 5\ 7\ \cdots \ 2i-1 \ h \ 0 \]
and $I(X)$ has a minimal free $R$-resolution of the following
type:

\[ 0\rightarrow R(-i-2)^{max(0,-\delta _{i+2})} \oplus  R(-i-3)^h
\rightarrow R(-i-1)^{max(0,\delta _{i+1})} \oplus
R(-i-2)^{max(0,\delta _{i+2})}
\]
\[\rightarrow R(-i)^{2i+1-h} \oplus
R(-i-1)^{max(0,-\delta _{i+1})}\oplus R(-2)\rightarrow R
\rightarrow R/I(X) \rightarrow 0
\]
where $\delta _n:=\Delta ^4 h_X(n)$ (the fourth difference of the
Hilbert function of
$X$).
\end{proposition}
\begin{proof} See \cite{GMR}; \S 4.
\end{proof}

\vskip 4mm
\begin{proposition}\label{res-rel-comp-gorODD}  Let $A$ be a general 
   Gorenstein Artinian graded algebra  of embedding dimension
$4$ and socle degree $2t+1$.  Assume that $A$ is relatively
compressed with respect to $\fa =(F)$ with $\deg(F)=2$ and $F$ general.
   If $2\le t$, then $A$ has a minimal free
   $R$-resolution of the following type:

\[
0 \rightarrow R(-2t-5) \rightarrow R(-t-4)^{2t+3}\oplus R(-2t-3)
\rightarrow R(-t-2)^{2t+3} \oplus  R(-t-3)^{2t+3}
\]
\[
\rightarrow R(-t-1)^{2t+3}\oplus R(-2) \rightarrow R \rightarrow A
\rightarrow 0.
\]

\end{proposition}
\begin{proof} We will explicitly construct a
   Gorenstein Artinian graded algebra  $A$ of embedding dimension
$4$, socle degree $2t+1$, relatively compressed with respect to
$\fa =(F)$, $\deg (F)=2$, and with the expected graded Betti
numbers. To this end, we consider a set of \linebreak $(t+1)^2+(t+1)$ general
points $X\subset\PP^3$ lying on  a smooth quadric $Q=V(F)\subset
\PP^3$. The $h$-vector of $X$ is
\[
1 \ 3 \ 5 \ \cdots \ 2t+1 \  \ t+1 \ 0 \
\]
    and $I(X)$ has a minimal free
$R$-resolution of the following type
\[
0 \rightarrow R(-t-4)^{t+1} \rightarrow
\begin{array}{c}
R(-t-3)^{t+2} \\
   \oplus \\
R(-t-2)^{t+1}
\end{array}
   \rightarrow
\begin{array}{c}
R(-t-1)^{t+2}\\
\oplus \\
R(-2)
\end{array} \rightarrow R
\rightarrow R/I(X) \rightarrow 0.
\]
   By Proposition \ref{MRC} such a set
of points $X$ exists.

By \cite{boij}; Theorem 3.2, the canonical module $\omega_X$ of
$R/I(X)$ can be embedded as an ideal $\omega_X (-2t-1)\subset
R/I(X)$ of initial degree $t+1$ and we have a short exact sequence
\[
0 \rightarrow \omega_X (-2t-1)\rightarrow R/I(X) \rightarrow A
\rightarrow 0
\]
where $A$ is a Gorenstein Artinian graded algebra of codimension
4, socle degree $2t+1$, $h$-vector
\[
1 \ 4 \ 9 \ 25 \ \cdots \ (t+1)^2 \ (t+1)^2 \ \cdots \ 25 \ 9 \ 4
\ 1.
\]
So, it is relatively compressed with respect to $\fa =(F)$,
$\deg(F)=2$. Moreover, applying the mapping cone process we get
that $A$ has the following minimal free $R$-resolution
\[
0 \rightarrow R(-2t-5) \rightarrow R(-t-4)^{2t+3}\oplus R(-2t-3)
\rightarrow R(-t-2)^{2t+3} \oplus  R(-t-3)^{2t+3}
\]
\[
\rightarrow R(-t-1)^{2t+3}\oplus R(-2) \rightarrow R \rightarrow A
\rightarrow 0.
\]
So it has the expected minimal free $R$-resolution in the sense
that the graded Betti numbers are the smallest consistent with the
Hilbert function of such relatively compressed Gorenstein graded
algebra.

\end{proof}

The approach used in Proposition \ref{res-rel-comp-gorODD} suggests that
there is a close relation between the Minimal Resolution Conjecture (MRC)
for points on a projective variety due to Musta\c t\v a (see \cite{MUS}
page 64) and the existence of relatively compressed Gorenstein algebras $A$
of odd socle degree and with the ``expected resolution'' in the sense that
the graded Betti numbers are the smallest consistent with the Hilbert
function of $A$.  We will end this section by writing down this relation.

\vskip 2mm \begin{definition}
Let $X\subset \PP^n$ be a projective variety of $\dim (X)\ge 1$. A
set $\Gamma $ of $\delta $ distinct points on $X$ is in {\em
general position} if
$$h_{R/I(\Gamma )}(t)= \min\{h_{R/I(X)}(t),\delta\}.$$
\end{definition}

If $X = \PP^n$ then the Minimal Resolution Conjecture predicts the graded Betti numbers of points in general position. It has been proved if the number of points is large compared to $n$ by Simpson and Hirschowitz \cite{Simps-Hirsch}, but may fail for a small number of points as shown by Eisenbud and Popescu \cite{Eisenb-P}.  If $X \neq \PP^n$ one has to modify the ``expectations.'' 
In \cite{MUS} pg.\ 64, Musta\c t\v a states the Minimal Resolution Conjecture
(MRC) for points on a projective variety. Let us  recall it.

\begin{conjecture} Let $X\subset \PP^n$ be a projective variety with
$d=\dim (X)\ge 1$, $\reg(X)=m$ and with Hilbert polynomial $P_X$.
Let $\delta $ be an integer with $P_X(r-1)\le \delta <P_X(r)$ for
some $r\ge m+1$ and let $\Gamma $ be a set of $\delta $  points on
$X$ in  general position. Let
$$
0 \rightarrow F_{n} \rightarrow F_{n-1}  \rightarrow \cdots
\rightarrow F_{2} \rightarrow F_{1}\rightarrow R \rightarrow R/I
\rightarrow 0
$$ be a minimal free $R$-resolution of $R/I(X)$. Then $R/I(\Gamma
)$ has a minimal free $R$-resolution of the following type

$$
0 \rightarrow F_{n}\oplus R(-r-n+1)^{a_{r+n-1}^n} \oplus
R(-r-n)^{a_{r+n}^n}\rightarrow F_{n-1} \oplus
R(-r-n+2)^{a_{r+n-2}^{n-1}} \oplus
R(-r-n+1)^{a_{r+n-1}^{n-1}}\rightarrow $$ $$ \cdots \rightarrow
F_{2} \oplus R(-r-1)^{a_{r+1}^2} \oplus
R(-r-2)^{a_{r+2}^2}\rightarrow F_{1}\oplus R(-r)^{a_{r}^1} \oplus
R(-r-1)^{a_{r+1}^1}\rightarrow R \rightarrow R/I(\Gamma )
\rightarrow 0
$$

\noindent with $a_{r+i}^{i}· a_{r+i}^{i+1}=0$ for $i=1, \cdots ,
n-1$.
\end{conjecture}

\begin{example} The MRC holds for $\delta  \gg 0$ points on a general
smooth rational quintic curve $C\subset \PP^3$ \cite{MUS}. 
\end{example}

For the purposes of this paper, it would be enough to know that the MRC
holds on a complete intersection variety.

\begin{proposition}\label{oddsocle}
    Let $A=R/I$ be a general 
     Gorenstein Artinian graded algebra  of embedding dimension
$n$ and socle degree $2t+1$.  Assume that $A$ is relatively
compressed with respect to a general complete intersection ideal
$\fa =(G_1, \cdots ,G_r)\subset k[x_1, \cdots , x_n]$, $r\le n-2$,
and $deg(G_1)=d_1, \cdots , deg(G_r)=d_r$. Set $m=\reg(X)$ where
$X=V(G_1, \cdots ,G_r)\subset \PP^{n-1}$ and $\underline{d}=(d_1,
\cdots , d_r)$. If $t\ge m$ and MRC holds for points on complete
intersection projective varieties, then $A$ has a minimal free
     $R$-resolution of the following type:

\[
0 \rightarrow R(-2t-1-n)\rightarrow F_{n-1} \rightarrow \cdots
\rightarrow F_2 \rightarrow F_1\rightarrow R \rightarrow A
\rightarrow 0 \]

\noindent where $F_i\cong F_{n-1}^{\vee }(-2t-1-n)$ for all $i=1,
\cdots , n-1$. Moreover, if $n$ is even (say, $n=2p$), then

\[
F_i= R(-t-i)^{\alpha _{i}(\underline{d},n,t)} \oplus
K_{i}(\underline{d}) \oplus K_{n-i}(\underline{d}) ^{\vee
}(-2t-1-n) \mbox{ for } i=1,\cdots , p-1; \mbox{ and }
\]

$$F_p=R(-t-p-1)^{\alpha _{p}(\underline{d},n,t)} \oplus
R(-t-p)^{\alpha _{p}(\underline{d},n,t)} \oplus
K_{p}(\underline{d}) \oplus K_{p}(\underline{d})^{\vee
}(-2t-1-n)$$

\vskip 2mm \noindent where $\alpha _i(\underline{d},n,t)$, $i=1,
\cdots p$, is completely determined by the Hilbert function of
$A$.

\vskip 2mm If $n$ is odd (say, $n=2p+1$), then

\[
F_i= R(-t-i)^{\alpha _{i}(\underline{d},n,t)} \oplus
K_{i}(\underline{d}) \oplus K_{n-i}(\underline{d})^{\vee
}(-2t-1-n) \mbox{ for } i=1,\cdots , p
\]

\vskip 2mm \noindent where $\alpha _i(\underline{d},n,t)$, $i=1,
\cdots p$, is completely determined by the Hilbert function of
$A$.

\end{proposition}

\begin{proof} We will only prove the case $n$ even and we leave to the
reader the case $n$ odd. To this end, we will explicitly construct
a   Gorenstein Artinian graded algebra  $A$ of embedding dimension
$n$, socle degree $2t+1$, relatively compressed with respect to
$\fa =(G_1, \cdots ,G_r)$ and with the expected graded Betti
numbers.

By hypothesis the  complete intersection projective variety
$X=V(G_1, \cdots ,G_r)\subset \PP^{n-1}$ satisfies  MRC  for any
set of $\delta $ points on $X$ in general position with
$P_X(s-1)\le \delta <P_X(s)$ for some $s\ge m+1$. Hence, since
$t\ge m$, for a suitable $\delta $ with $P_X(t)\le \delta
<P_X(t+1)$, a set $\Gamma $ of $\delta $ general points on $X$ has
a minimal free $R$-resolution of the following type:

\begin{equation} \label{resol1}0 \rightarrow
K_{n-1}(\underline{d}) \oplus R(-t-n)^{a_{t+n}^{n-1}}\rightarrow
K_{n-2}(\underline{d}) \oplus R(-t-n+1)^{a_{t+n-1}^{n-2}}
\end{equation}
$$ \rightarrow \cdots \rightarrow K_{n/2}(\underline{d})
\oplus R(-t-\frac{n}{2})^{a_{t+\frac{n}{2}}^{n/2}} \oplus
R(-t-\frac{n}{2}-1)^{a_{t+\frac{n}{2}+1}^{n/2}} \rightarrow \cdots
  $$
  $$  \rightarrow
K_{2}(\underline{d}) \oplus R(-t-2)^{a_{t+2}^2} \oplus \rightarrow
K_{1}(\underline{d})\oplus  R(-t-1)^{a_{t+1}^1}\rightarrow R
\rightarrow R/I(\Gamma ) \rightarrow 0.
$$

\vskip 2mm Dualizing and twisting the  exact sequence
(\ref{resol1}), we get a minimal free $R$-resolution of the
canonical module $\omega _{\Gamma }(-2t-1)$ of $\Gamma $: \vskip
2mm
\begin{equation} \label{resol2}
\begin{array}{c}
0 \rightarrow R(-2t-1-n) \rightarrow R(-t-n)^{a _{t+1}^1} \oplus
K_{1}(\underline{d})^{\vee }(-2t-1-n)\rightarrow \\ \vbox{\vskip .8cm}
  R(-t-n+1)^{a_{t+2}^2 }\oplus K_{2}(\underline{d}^{\vee }(-2t-1-n)
\rightarrow \cdots \rightarrow\\ \vbox{\vskip .8cm}
  R(-t-1-\frac{n}{2})^{a_{t+\frac{n}{2}}^{n/2}} \oplus R(-t-\frac{n}{2})^{a
_{t+\frac{n}{2}+1}^{n/2}} \oplus K_{n/2}(\underline{d})^{\vee
}(-2t-1-n) \rightarrow \cdots \rightarrow \\ \vbox{\vskip .8cm}
    R(-t-2)^{a_{t+n-1}^{n-2}} \oplus K_{n-2}(\underline{d})^{\vee
}(-2t-1-n)\rightarrow \\ \vbox{\vskip .8cm}
  R(-t-1)^{a _{t+n}^{n-1}} \oplus
K_{n-1}(\underline{d})^{\vee }(-2t-1-n)\rightarrow \omega _{\Gamma
}(-2t-1) \rightarrow 0.
\end{array}
\end{equation}

By \cite{boij}; Theorem 3.2, the canonical module $\omega_{\Gamma
}$ of $R/I(\Gamma )$ can be embedded as an ideal $\omega_{\Gamma }
(-2t-1)\subset R/I(\Gamma )$ of initial degree $t+1$ and we have a
short exact sequence
\begin{equation} \label{res3}
0 \rightarrow \omega_{\Gamma } (-2t-1)\rightarrow R/I(\Gamma )
\rightarrow A \rightarrow 0
\end{equation}
where $A$ is a Gorenstein Artinian graded algebra of codimension
$n$, socle degree $2t+1$.  A straightforward computation using the exact
sequences (\ref{resol1}), (\ref{resol2}) and (\ref{res3}) gives us
$h_A(\lambda )= \min\{h_{R/\fa}(\lambda ), h_{R/\fa}(2t+1-\lambda )\}$.
Therefore, $A$ is relatively compressed with respect to
$\fa =(G_1, \cdots ,G_r)$. Moreover, using the exact sequences
(\ref{resol1}), (\ref{resol2}) and (\ref{res3}) and applying the
mapping cone process, we get that $A$ has the following minimal
free $R$-resolution:

\[
0 \rightarrow R(-2t-1-n)  \rightarrow
R(-t-n)^{a^1_{t+1}+a^{n-1}_{t+n}}\oplus K_{n-1}(\underline{d})
\oplus K_{1}(\underline{d})^{\vee }(-2t-1-n) \rightarrow
\]
\[R(-t-n-1)^{a^2_{t+2}+a^{n-2}_{t+n-1}}\oplus K_{n-2}(\underline{d})
\oplus K_{2}(\underline{d})^{\vee }(-2t-1-n) \rightarrow
  \cdots \rightarrow
\]
\[R(-t-\frac{n}{2})^{a^{n/2}_{t+\frac{n}{2}}+a^{n/2}_{t+\frac{n}{2}+1}}\oplus
R(-t-\frac{n}{2}-1)^{a^{n/2}_{t+\frac{n}{2}}+a^{n/2}_{t+\frac{n}{2}+1}}\oplus
  K_{n/2}(\underline{d})
\oplus K_{n/2}(\underline{d})^{\vee }(-2t-1-n)
\]
\[
\rightarrow
  \cdots \rightarrow
R(-t-2)^{a^2_{t+2}+a^{n-2}_{t+n-1}}\oplus K_{2}(\underline{d})
\oplus K_{n-2}(\underline{d})^{\vee }(-2t-1-n) \rightarrow
\]
\[R(-t-1)^{a^1_{t+1}+a^{n-1}_{t+n}}\oplus K_{1}(\underline{d})
\oplus K_{n-1}(\underline{d})^{\vee }(-2t-1-n) \rightarrow
  R \rightarrow A \rightarrow 0.
\]

\vskip 2mm So $A$ has the expected minimal free $R$-resolution in
the sense that the graded Betti numbers are the smallest
consistent with the Hilbert function of such relatively compressed
Gorenstein algebra and there are no non-Koszul ghost terms.
\end{proof}


\section{Unexpected behavior for  the minimal free resolution of relatively
compressed level algebras of socle dimension $c \geq 2$: ``ghost'' terms}
\label{sec-Conj}

We saw in section \ref{sec-exist} that the Hilbert function of a relatively
compressed algebra can fail to be the ``expected" one.  In this
section we consider
the question of when (if at all) the minimal free resolution of a relatively
compressed level algebra must have redundant (or ``ghost'') terms that are not
computable only from the Hilbert function.

We first give two examples of level $c\ge 2$
    algebras relatively compressed with
respect to a complete intersection $\fc $,  with fixed socle
degree $s$, embedding dimension 3 and with the ``expected" minimal
free $R$-resolution.

\vskip 2mm The first example involves a {\tt macaulay} computation, while
the second
   does not.

\begin{example} We set $R=k[x,y,z]$.
Let $A$ be a level 2 graded algebra of embedding dimension 3,
socle degree 15 and relatively compressed with respect to a
complete intersection ideal  $\fa =(F_1,F_2,F_3)$, $\deg(F_{i})=9$
for $i=1,2,3$. Note that the Hilbert function of $R/\fa$ is
\[
1 \ 3 \ 6 \ 10 \ 15 \ 21 \ 28 \ 36 \ 45 \ 52 \ 57 \ 60 \ 61 \ 60 \ 57
\ 52 \ 45 \ 36
\ 28 \ 21 \ 15 \ 10 \ 6 \ 3 \ 1.
\]
The expected $h$-vector of $A$ (according to the discussion
preceding Example \ref{orig 2.4}) is
\begin{equation}\label{predicted1}
1 \  3 \  6 \ 10 \ 15 \ 21 \ 28 \ 36
\ 45 \ 52\ 42 \ 30 \ 20 \ 12
\ 6 \ 2;
\end{equation}
    and the expected minimal free $R$-resolution is

\[
0 \rightarrow R(-18)^{2} \rightarrow R(-11)^{15}\oplus R(-12)^{4}
\rightarrow R(-9)^{3} \oplus R(-10)^{15} \rightarrow R \rightarrow
A \rightarrow 0.
\]

Let us explicitly construct it. To this end, we consider the ideal
$J=\fa +(F_4,F_5)$ where  $F_4$ and $F_5$ are two general forms of
degree 9.  We know the Hilbert function of $R/J$ thanks to
\cite{anick} (using a
calculation similar to that in Example \ref{orig 2.6}):
\[
1 \ 3 \ 6 \ 10 \ 15 \ 21 \ 28 \ 36 \ 45 \ 50 \ 51 \ 48 \ 41 \ 30 \ 15.
\]
Using {\tt macaulay} we compute a minimal free $R$-resolution
of $J$ and we get

\[
0 \rightarrow R(-17)^{15} \rightarrow R(-16)^{15}\oplus R(-15)^{4}
\rightarrow R(-9)^{5} \rightarrow R \rightarrow R/J \rightarrow 0.
\]

The ideal $\fa $ links $J$ to an ideal $J':=[\fa : J]$ and it is
easily seen that $A=R/J'$ is a level 2 graded algebra of embedding
dimension 3 and socle degree 15.  One quickly checks that its Hilbert
function is
the one predicted above in (\ref{predicted1}), so $A$ is relatively
compressed with
respect to a complete intersection ideal  $\fa =(F_1,F_2,F_3)$,
$deg(F_{i})=9$ for
$i=1,2,3$.
Using the standard mapping cone construction we get that $A$ has
the following minimal free $R$-resolution:
\[
0 \rightarrow R(-18)^{2} \rightarrow R(-11)^{15}\oplus R(-12)^{4}
\rightarrow R(-9)^{3} \oplus R(-10)^{15} \rightarrow R \rightarrow
A \rightarrow 0.
\]
So, it has the expected minimal free $R$-resolution in the sense
that the graded Betti numbers are the smallest consistent with the
Hilbert function of such relatively compressed level  graded
algebras of socle dimension 2. Note that the only place where we needed the aid
of the computer
was to confirm the expected minimal free resolution of an ideal of five general
forms of degree 9.
\end{example}

\begin{example} We set $R=k[x,y,z]$.
Let $A$ be a level  graded algebra of socle dimension  4 and embedding dimension 3,
socle degree 5 and relatively compressed with respect to an ideal
$\fa =(F)$, $\deg(F)=2$. The $h$-vector of $A$ is 1 3 5 7 9 4; and
the expected minimal free $R$-resolution is

\[
0 \rightarrow R(-8)^{4} \rightarrow R(-6)^{8}\oplus R(-7)^{3}
\rightarrow R(-2) \oplus R(-5)^{7} \rightarrow R \rightarrow A
\rightarrow 0.
\]

Let us explicitly construct it. To this end, we consider a set of
29 general points $X\subset \PP^3$ lying on a smooth quadric
surface $Q\subset \PP^3$. The $h$-vector of $X$ is 1 3 5 7 9 4;
and $I(X)\subset S:=k[x_0,x_1,x_2,x_3]$ has a minimal free
$S$-resolution of the following type (cf. \cite{GMR})

\[
0 \rightarrow S(-8)^{4} \rightarrow S(-6)^{8}\oplus S(-7)^{3}
\rightarrow S(-2) \oplus S(-5)^{7} \rightarrow S \rightarrow
S/I(X) \rightarrow 0.
\]

The Artinian reduction of $S/I(X)$ is a level  graded Artinian
algebra $A$ of socle dimension 4 and embedding dimension 3, socle degree 5 and
relatively compressed with respect to an ideal $\fa =(F)\subset
k[x,y,z]$, $\deg(F)=2$. The $h$-vector of $A$ is 1 3 5 7 9 4; and
the  minimal free $R$-resolution of $A$ is

\[
0 \rightarrow R(-8)^{4} \rightarrow R(-6)^{8}\oplus R(-7)^{3}
\rightarrow R(-2) \oplus R(-5)^{7} \rightarrow R \rightarrow A
\rightarrow 0.
\]
Again, it has the expected minimal free $R$-resolution in the
sense that the graded Betti numbers are the smallest consistent
with the Hilbert function of such relatively compressed level 
graded algebra with socle dimension 2.
\end{example}

   In the next example we show that sometimes there must be redundant
terms in the
minimal free resolution of a relatively compressed level algebra.

\begin{example}\label{another ghost}
We work in four variables, $R = k[x_1,x_2,x_3,x_4]$.  Let $I_1$ be a
general ideal
with $h$-vector 1 2, so its minimal free resolution is
\[
0 \rightarrow R(-5)^2 \rightarrow R(-4)^7 \rightarrow R(-3)^8 \oplus R(-2)
\rightarrow R(-1)^2 \oplus R(-2)^3 \rightarrow I_1 \rightarrow 0.
\]
Linking with a complete intersection of type $(3,3,3,3)$ we get an
ideal $I_2$ with
$h$-vector
\[
1 \ 4 \ 10 \ 16 \ 19 \ 16 \ 10 \ 2
\]
   and minimal free resolution
\[
0 \rightarrow
\begin{array}{c}
R(-11)^2 \\
\oplus \\
   R(-10)^3
\end{array}
\rightarrow
\begin{array}{c}
R(-9)^{12} \\
\oplus \\
R(-10)
\end{array}
\rightarrow
\begin{array}{c}
R(-8)^7 \\
\oplus \\
R(-6)^6
\end{array}
\rightarrow
\begin{array}{c}
R(-7)^2 \\
\oplus \\
R(-3)^4
\end{array}
\rightarrow I_2 \rightarrow 0
\]
(that you can compute from that of $I_1$ with the mapping cone).
$I_2$ is an ideal of general forms of degrees 3,3,3,3,7,7, and its minimal free
resolution has a non-Koszul ghost term $R(-10)$.  This approach for
finding ideals
of general forms with non-Koszul ghost terms was developed in \cite{MMR3}.

Now we link $I_2$ with a complete intersection $\fc$ of type
$(3,3,7,7)$.  The Hilbert function of $\fc$ is
\[
1 \ 4 \ 10 \ 18 \ 27 \ 36 \ 45 \ 52 \ 55 \ 52 \ 45 \ 36 \ 27 \ 18 \ 10 \ 4 \ 1.
\]
Letting $I_3$ be the residual, its Hilbert function is
\[
1 \ 4 \ 10 \ 18 \ 27 \ 36 \ 45 \ 52 \ 55 \ 50 \ 35 \ 20 \ 8 \ 2.
\]
Note first that $I_3$ does not ``look'' relatively compressed in $\fc$:
\[
\begin{array}{lccccccccccccccccccccccccccccccccccc}
\deg & 0 & 1 & 2 & 3 & 4 & 5 & 6 & 7 & 8 & 9 & 10 & 11 & 12 & 13 & 14
& 15 & 16 \\
\fc & 1 & 4 & 10 & 18 & 27 & 36 & 45 & 52 & 55 & 52 & 45 & 36 & 27 &
18 & 10 & 4 & 1
\\
I_3 & 1 & 4 & 10 & 18 & 27 & 36 & 45 & 52 & 55 & 50 & 35 & 20 & 8 & 2
\end{array}
\]
   The 35 ``should'' be a 36, and the 50 ``should'' be 52.  But given
any relatively
compressed level algebra of socle dimension  2 and socle degree 13 inside a
complete intersection
of type
$(3,3,7,7)$, the residual is an ideal generated by six forms, of degrees
3,3,3,3,7,7.  The smallest Hilbert function of such an ideal is given
by our $I_2$
(general forms), so the biggest (i.e.\ relatively compressed) for the
level algebra
is the one above.

But now we consider minimal free resolutions.  First note that in the
resolution of
$I_2$ above, all the copies of $R(-6)$ are Koszul.  We know the
resolution of $\fc$
(the Koszul resolution).  The mapping cone gives the following.  Note
that we split
off not only generators of degrees 3,3,7,7, but also one first syzygy
$R(-6)$, namely
the Koszul one.  Studying the mapping cone carefully, we see that there
is no other possible splitting.  We get:
\[
0 \rightarrow R(-17)^2 \rightarrow
\begin{array}{c}
R(-12)^7 \\
\oplus \\
R(-14)^5
\end{array}
\rightarrow
\begin{array}{c}
R(-11)^{12} \\
\oplus \\
R(-10)^5 \\
\oplus \\
R(-6)
\end{array}
\rightarrow
\begin{array}{c}
R(-10)^3 \\
\oplus \\
R(-9)^2 \\
\oplus \\
R(-7)^2 \\
\oplus \\
R(-3)^2
\end{array}
\rightarrow I_3 \rightarrow 0.
\]
Note that four of the copies of $R(-10)$ in the second free module
are Koszul, but
the fifth one is not.  So we have a relatively compressed
level algebra with a non-Koszul ghost term at the beginning of the resolution.
\end{example}

The following result generalizes the last example.  It is based on \cite{MMR3}
Theorem 3.3.

\begin{proposition} \label{link ghost}
Let $R = k[x_1,\dots,x_n]$ and let $\fc' = (F_1,\dots,F_n)$ be a complete
intersection of forms of degree $d_1,\dots,d_n$, respectively.  We do
not assume
that $d_1 \leq \dots \leq d_n$, but we do assume that each $d_i > 2$.
Suppose that
$d_1 = \dots = d_p = a$ (say) for some $1 \leq p \leq n-2$.  Let $d =
d_1 + \dots
+d_n$ and choose general forms $F_{n+1},\dots,F_{n+p}$ all of degree
$d-n-1$.  Let
$\fc$ be the complete intersection
$(F_{p+1},\dots,F_n,F_{n+1},\dots,F_{n+p})$ and
let $J = (F_1,\dots,F_{n+p})$.  Let $e = \sum_{i=p+1}^n d_i + p(d-n-1)$.
Then the residual ideal $I = \fc : J$ is a level
algebra of socle dimension $p$ and socle degree $e-n-a$, and
$R/I$ is relatively compressed in $R/\fc$.  Furthermore, the minimal
free resolution
of $R/I$
\[
0 \rightarrow \mathbb F_n \rightarrow \dots \rightarrow \mathbb F_2 \rightarrow
\mathbb F_1 \rightarrow R \rightarrow R/I \rightarrow 0
\]
has ghost terms between $\mathbb F_j$ and $\mathbb F_{j+1}$ for $1 \leq j \leq
n-p-1$.
\end{proposition}

\begin{proof}
Clearly $R/J$ has a minimal free resolution
\[
0 \rightarrow \mathbb G_n \rightarrow \dots \rightarrow \mathbb
G_{p+2} \rightarrow
\mathbb G_{p+1} \rightarrow \dots \rightarrow \mathbb G_1 \rightarrow
R \rightarrow
R/J \rightarrow 0
\]
where $\mathbb G_1 = R(-a)^p \oplus \bigoplus_{i=p+1}^n R(-d_i)
\oplus R(-d+n+1)^p$.
Furthermore, by  \cite{MMR3} Theorem 3.3, for $j = p+1,\dots,n-1$
there is a ghost
term $R(-d+n+1-j)$ between $\mathbb G_j$ and $\mathbb G_{j+1}$ that
does not arise
from Koszul syzygies.  Note also that since the socle degree of $R/J$
is $d-n-1$,
the largest twist of any $\mathbb G_j$ (including $\mathbb G_n$) is
\[
\begin{array}{rcl}
R(-d+1-n+j) & =
   & R(-d_1-\dots-d_n +1-n+j) \\
& = & R(-pa - d_{p+1} - \dots -d_n +1-n+j).
\end{array}
\]

Consider the minimal free (Koszul) resolution of $\fc$
\[
\begin{array}{cccccccccccccccccccccc}
&&&&&&&&&& 0 \\
&&&&&&&&&& \downarrow \\
0 & \rightarrow & \mathbb K_n & \rightarrow  \dots  \rightarrow &
\mathbb K_{p+2}
& \rightarrow & \mathbb K_{p+1} & \rightarrow \dots \rightarrow & \mathbb K_1 &
\rightarrow & \fc & \rightarrow & 0 \\
&& \downarrow && \downarrow && \downarrow && \downarrow && \downarrow \\
0 & \rightarrow & \mathbb G_n & \rightarrow  \dots  \rightarrow &
\mathbb G_{p+2}
& \rightarrow & \mathbb G_{p+1} & \rightarrow \dots \rightarrow & \mathbb G_1 &
\rightarrow & J & \rightarrow & 0
\end{array}
\]
Using the mapping cone procedure, we
obtain a free resolution for $R/I$.  We include the obvious splitting
of the terms
$R(-a)^p$ between $\mathbb K_1$ and $\mathbb G_1$.
\[
0 \rightarrow R(-e+a)^p \rightarrow \mathbb G_2^\vee (-e) \rightarrow \dots
\rightarrow
\mathbb G_n^\vee (-e) \oplus \mathbb K_{n-1}^\vee \rightarrow R \rightarrow R/I
\rightarrow 0.
\]
It is clear that no further terms split from $R(-e+a)^p$ (by the
minimality of the
resolution of $J$).  Since the ghost terms of the minimal free
resolution of $R/J$
are not Koszul and since the generators of $\fc$ are taken from the
generators of
$J$, it is clear that no splitting will remove these terms after
taking the mapping
cone, so they remain (suitably twisted) in the minimal free
resolution of $R/I$.
\end{proof}

\begin{example}
In Proposition \ref{link ghost} we concluded that under certain hypotheses
we can find a relatively compressed level algebra with non-Koszul ghost terms
at the beginning of the resolution.  Following \cite{MMR3} Corollary 3.13, we
can even arrange some splitting at the beginning of the resolution, leaving
ghost terms only in the middle.  Note that (according to the observation
following \cite{MMR3} Corollary 3.13) this will only work for $n=4,5,6$.
Here is one example.

Let $n=5$ and start with a quotient of $R$ with Hilbert function $1 \ 1$
(it is a complete intersection of type $(1,1,1,1,2)$).  We link using a
complete intersection of type $(2,4,4,4,4)$ to obtain an ideal with
generators of degrees 2,4,4,4,4,12, with ``expected" Hilbert function.  We
link this in turn using generators of degrees 2,4,4,4,12 to obtain a
Gorenstein ideal $J$ with Hilbert function
\[
1 \ \ 5 \ \ 14 \ \ 30 \ \ 52 \ \ 76 \ \ 98 \ \ 114 \ \ 123 \ \ 123 \ \
114 \
\ 98 
\  76 \
\ 52 \
\ 30 \ \ 14 \ \ 5 \ \ 1.
\]
The minimal free resolution of $J$ can be computed (although
it is very tedious).  It is
\[
0 \rightarrow R(-22) \rightarrow
\begin{array}{c}
R(-20) \\
\oplus \\
R(-18)^3 \\
\oplus \\
R(-13)^4
\end{array}
\rightarrow
\begin{array}{c}
R(-16)^3 \\
\oplus \\
R(-14)^3 \\
\oplus \\
R(-12)^6 \\
\oplus \\
R(-11)^{4} \\
\oplus \\
R(-10)^3
\end{array}
\rightarrow
\begin{array}{c}
R(-12)^3 \\
\oplus\\
R(-11)^4 \\
\oplus \\
R(-10)^6 \\
\oplus \\
R(-8)^3 \\
\oplus \\
R(-6)^3
\end{array}
\rightarrow
\begin{array}{c}
R(-9)^4 \\
\oplus \\
R(-4)^3 \\
\oplus \\
R(-2)
\end{array}
\rightarrow J \rightarrow 0. 
\]
Then we point out first that there are no ghost terms, Koszul or otherwise,
between the first free module and the second.  Furthermore, there are some
copies of $R(-10)$ between the second and third modules, but the three
copies of $R(-10)$ in the third module are Koszul.  But now consider the
copies of $R(-12)$ between the second and third modules.  None of the copies
of $R(-12)$ in the second module come because of Koszul syzygies, and in the
third module at most one copy of $R(-12)$ is Koszul.  So even after
accounting for that, we have non-Koszul ghost terms between the second and
third modules in the resolution.
\end{example}

\begin{example} \label{gor ex}
When $n=3$ we are not aware of any relatively compressed level algebras of
socle dimension $c \geq 2$ that have non-Koszul ghost terms.  However, for
$c=1$ these do exist.  We heavily use the results of \cite{MMR2}.
First note that if $I = (G_1,G_2,G_3,G_4)$ is an ideal of general forms in 3
variables, and if $\fc$ is a complete intersection defined by any three of
the four generators, then the residual $\fc : I = G$ is a Gorenstein ideal.
Furthermore, we claim that it is relatively compressed with respect to
$\fc$.  Indeed, any Gorenstein ideal containing a complete intersection of
the same degrees, and having the same socle degree, will be linked by that
complete intersection to an ideal of four forms of the same degrees as
$G_1,G_2,G_3,G_4$.  Since
$R/I$ has the smallest possible Hilbert function among all such ideals, by
linkage
$R/G$ has the largest Hilbert function among all Gorenstein ideals with that
socle degree, containing a complete intersection of that type.

In the case where $\deg G_1 \leq \dots \leq \deg G_4$ and $\fc =
(G_1,G_2,G_3)$, it was shown in
\cite{MMR2} Proposition 4.1 that there appears a non-Koszul ghost term if
and only if either (i) $d_2+d_3< d_1 +d_4+3$ and $d_1+d_2+d_3+d_4$ is even,
or (ii) $d_2+d_3 \geq d_1 +d_4+3$ and $d_1+d_2+d_3+d_4$ is even.  In
particular, such examples exist.
\end{example}

We believe that ``most of the time," a relatively compressed level
quotient of a complete intersection has only Koszul ghost terms.  We believe,
furthermore, that the only counterexamples arise through linkage, as
special cases
of Proposition \ref{link ghost}.  More precise conjectures are the following:

\begin{conjecture}
Let $\fc$ be a general complete intersection of fixed generator degrees in a
polynomial ring
$R =  k[x_1,\dots,x_n]$, and
let $A = R/I$ be a relatively compressed level quotient of $R/\fc$ of
socle degree
$\geq 2$ and socle dimension  $c \geq 1$.  Assume that either $n=3$
and $c\geq 2$ 
or else $n \geq 4$.  If the minimal free resolution of
$R/I$ has non-Koszul  ghost terms then
$I$ is linked in two steps, first  by $\fc$ and then by a
``predictable'' complete
intersection, to an ideal containing at least two independent linear forms.
\end{conjecture}

Note that the ideals in Example \ref{another ghost} and Proposition
\ref{link ghost} have this property of being linked in two steps to an ideal
containing at least two independent linear forms.  Indeed, the Koszul
relations of this latter ideal are what produce the non-Koszul relations of
the final ideal.

\begin{conjecture} \label{2nd conj}
Assume that $\fc$ is a complete intersection generated by forms all of the
same degree.  Let $A$ be a relatively compressed level
Artinian quotient of
$R/\fc$ of socle dimension  $c$ and socle degree $\geq 2$,and assume
that either $n=3$ 
and $c \geq 2$, or $n \geq 4$.  Then the minimal free resolution of $A$  has
no ghost terms (Koszul or otherwise).
\end{conjecture}

\begin{remark}
The assumptions in Conjecture \ref{2nd conj} are necessary.  In the case of
socle degree 1, a counterexample would be any algebra with Hilbert function
$1 \ t \ 0$ for $t \leq n-2$.  In the case of Gorenstein algebras of height
3, we have from \cite{MMR2} Proposition 4.1 that four general forms all of
the same degree are {\em always} linked (using three of the four generators)
to a Gorenstein ideal whose minimal free resolution has a non-Koszul ghost
term.  That this ideal is relatively compressed follows for instance from
\cite{MMR2} Lemma 2.6.
\end{remark}

\begin{remark} \label{no koszul}
Let $R = k[x_1,\dots,x_n]$ be a polynomial ring where the Fr\"oberg
conjecture holds (equivalently, where Maximal Rank Property holds -- cf.\
\cite{MMR3}).  For instance, this is true for $n=2,3$ and conjecturally for
all $n$.  Let $I$ be an ideal minimally generated by $r \geq n+1$ general
forms of degrees $a_1 \leq a_2 \dots \leq a_r$.  In particular, we are
assuming that none of these generators is redundant.  Let $J$ be the complete
intersection defined by the first three generators.  Then $R/I$ is a
quotient of $R/J$, and in particular (since $h_{R/J}(a_1+a_2+\dots+a_n -n) =
1$), the socle degree, $\delta$, of $R/I$ is strictly less than that of
$R/J$.  That is,
\[
\delta +n < a_1+a_2+\dots +a_n.
\]
It follows that if the minimal free resolution of $R/I$ is
\[
0 \rightarrow F_n \rightarrow F_{n-1} \rightarrow \dots \rightarrow F_1
\rightarrow R \rightarrow R/I \rightarrow 0
\]
then no summand of $F_n$ corresponds to a Koszul $(n-1)$-st syzygy among any
of the generators, since the highest twist of $F_n$ is $R(-\delta-n)$.
\end{remark}

\begin{proposition} \label{resol in n vars}
Let $R = k[x_1,\dots,x_n]$ and assume that $R$ satisfies the Fr\"oberg
conjecture (i.e.\ any ideal of general forms  satisfies the
Maximal Rank Property).  Let
$A$ be a level quotient of
$R$ of socle dimension 
$c\ge 2$. (The case
$c=1$ is the Gorenstein case that we have already discussed in Example
\ref{gor ex}).  Assume that $A$ has socle degree $s$ and is relatively
compressed with respect to a general complete intersection $\fc
=(F_1,F_2, \dots, F_n)\subset R$. Set $d_i=deg(F_i)$,  and
$d=d_1+d_2+ \dots +d_n$.  Let $I$ be the ideal generated by $\fc$ and $c$
general forms of degree $d-s-n$, and set
$\delta$ to be the socle degree of $R/I$.  Let $\mathbb G = R(-d_1) \oplus
R(-d_2) \oplus \dots \oplus R(-d_n) \oplus R(s+3-d)^c$ and let $K_t =
\wedge^t
\mathbb G$.   Finally, let $\mathbb
F = R(-d_1) \oplus \dots \oplus R(-d_n)$ and let $L_t = \wedge^t \mathbb F$.

Then $A$ has a
free $R$-resolution  of the following type:  {\footnotesize
\[
0 \rightarrow R(-s-n)^c  \rightarrow
\begin{array}{c}
  R(\delta +1-d)^{z_2} \\ \oplus \\ R(\delta
+2-d)^{w_2} \\
\oplus \\ ((K_2)^{\le \delta})^{\vee }(-d)
\end{array}
\rightarrow
\begin{array}{c}
R(\delta+2-d)^{z_3} \\
\oplus \\
R(\delta+3-d)^{w_3} \\
\oplus \\
((K_3)^{\leq \delta+1})^\vee (-d)\\
\oplus \\
(L_2)^\vee (-d)
\end{array}
\rightarrow
\dots \rightarrow \hskip -.4cm
\begin{array}{c}
R(\delta +n-1-d)^{z_n} \\ \oplus \\ R(\delta +n-d)^{w_n} \\ \oplus \\
R(-d_1) \\ \oplus \\ \vdots \\ \oplus \\ R(-d_n)
\end{array}
\rightarrow R \rightarrow A \rightarrow 0
\]
}
where $w_n$ and $z_2$ are determined by the Hilbert Function of $A$.

\end{proposition}

\begin{proof}  We have
$\fc=(F_1,F_2,\dots, F_n)$,  $I=(\fc ,G_1,\cdots ,G_c)$ with
$\deg(G_j)=d-s-n$ for  all $j=1, \cdots , c$.  Note that the socle degree
$\delta $ is determined by $d_1,d_2, \dots, d_n, c$ and
$d-s-n$ because we have assumed that an ideal of general forms in $R$
satisfies the Maximal Rank Property.)

By \cite{MMR3} Theorem 3.15 and Remark \ref{no koszul} above, $I$ has a  free
$R$-resolution  of the following type:
{\footnotesize
\begin{equation} \label{resol of gen forms}
0 \rightarrow
\begin{array}{c}
R(-\delta -n+1)^{z_n} \\ \oplus \\ R(-\delta -n)^{w_n}
\end{array} \rightarrow
\begin{array}{c}
R(-\delta-n+2)^{z_{n-1}} \\ \oplus \\ R(-\delta-n+1)^{w_{n-1}} \\
\oplus \\ (K_{n-1})^{\le \delta+n-3}
\end{array}
\rightarrow
\dots \rightarrow
\begin{array}{c}
R(-\delta-1)^{z_2} \\
\oplus \\
R(-\delta-2)^{w_2} \\
\oplus \\
(K_2)^{\leq \delta}
\end{array}
\rightarrow
\begin{array}{c}
R(s+n-d)^c \\ \oplus \\ R(-d_1) \\ \oplus \\ \vdots \\ \oplus \\ R(-d_n)
\end{array}
\rightarrow I \rightarrow 0.
\end{equation}
}

Note that $\fc$ has a minimal free (Koszul) resolution
\[
0 \rightarrow L_n \rightarrow L_{n-1} \rightarrow \dots \rightarrow L_1
\rightarrow \fc \rightarrow 0
\]
where $L_n = R(-d)$.
Applying the mapping cone process we get that $J=[\fc :I]$
has a  free $R$-resolution  of the following type: {\footnotesize
\[
0 \rightarrow R(-s-n)^c  \rightarrow
\begin{array}{c}
  R(\delta +1-d)^{z_2} \\ \oplus \\ R(\delta
+2-d)^{w_2} \\
\oplus \\ ((K_2)^{\le \delta})^{\vee }(-d)
\end{array}
\rightarrow
\begin{array}{c}
R(\delta+2-d)^{z_3} \\
\oplus \\
R(\delta+3-d)^{w_3} \\
\oplus \\
((K_3)^{\leq \delta+1})^\vee (-d)\\
\oplus \\
(L_2)^\vee (-d)
\end{array}
\rightarrow
\dots \rightarrow
\begin{array}{c}
R(\delta +n-1-d)^{z_n} \\ \oplus \\ R(\delta +n-d)^{w_n} \\ \oplus \\
R(-d_1) \\ \oplus \\ \vdots \\ \oplus \\ R(-d_n)
\end{array}
\rightarrow J \rightarrow 0
\]
}
as claimed.
\end{proof}

\begin{remark}
It is natural to wonder how close the resolution in Proposition \ref{resol
in n vars} is to being minimal.  The first consideration  is whether
(\ref{resol of gen forms}) is minimal.  The minimality of the first free
module is completely determined by the Maximal Rank Property, since the
forms are general.  When the redundant generators are removed, it is then
possible (in any given example) to determine how much splitting occurs in the
mapping cone.  So in fact, the only unknowns concern the values of the
graded Betti numbers in redundant terms.  We conjecture that when $n=3$, one
or the other must always be zero (i.e. $yz   = 0$ in the following result),
so there are no non-Koszul ghost terms.  However, for $n \geq 4$ we have
seen in Example \ref{another ghost} that this is not true.
\end{remark}

\begin{corollary} \label{resol in 3 vars}
Let $R = k[x,y,z]$ and let $A$ be a level quotient of $R$ of socle
dimension  $c\ge 2$.
(The case
$c=1$ is the Gorenstein case that we have already discussed in Example
\ref{gor ex}).  Assume that $A$ has socle degree $s$ and is relatively
compressed with respect to a general complete intersection $\fc
=(F_1,F_2,F_3)\subset R$. Set $d_i=deg(F_i)$,  and
$d=d_1+d_2+d_3$.  Let $I$ be the ideal generated by $\fc$ and $c$ general
forms of degree $d-s-3$, and set
$\delta$ to be the socle degree of $R/I$.  Let $\mathbb F = R(-d_1) \oplus
R(-d_2) \oplus R(-d_3) \oplus R(s+3-d)^c$ and let $K_2 = \wedge^2 \mathbb
F$.

Then $A$ has a
free $R$-resolution  of the following type:
\[
0 \rightarrow R(-s-3)^c  \rightarrow
\begin{array}{c}
R(\delta +1-d)^x \\ \oplus \\ R(\delta
+2-d)^y \\
\oplus \\ (K_2)^{\le \delta})^{\vee }(-d)
\end{array}
\rightarrow
\begin{array}{c}
  R(\delta +2-d)^z \\ \oplus \\ R(\delta +3-d)^w \\ \oplus \\
R(-d_1) \\ \oplus \\ R(-d_2) \\ \oplus \\ R(-d_3)
\end{array}
\rightarrow R \rightarrow A \rightarrow 0
\]
where $w$ and $x$ are determined by the Hilbert Function of $A$.
\end{corollary}

\begin{proof}
This follows from Proposition \ref{resol in n vars} since $k[x,y,z]$
satisfies Fr\"oberg's conjecture (cf.~\cite{anick}).
\end{proof}


\section{Applications to ideals of general forms}

\label{sec-Applications}

In this section, as an application of Theorem
\ref{res-rel-comp-gorcodim-n}, we get new results about the
generic graded Betti numbers of an almost complete intersection
ideal. The idea (which is not new) is to link an Artinian
Gorenstein graded algebra $R/G$, relatively compressed with respect
to a complete intersection ideal $\fa =(G_1, \cdots ,G_r)\subset
k[x_1,\cdots ,x_n]$, to an almost complete intersection ideal $(\fa
,G_{r+1},\cdots ,G_n,G_{n+1 })$ via a complete intersection $(\fa
,G_{r+1},\cdots ,G_n)$ where $G_{r+1}$, $\cdots $, $G_n$ are
suitably chosen.

\vskip 2mm Let $R = k[x_1,\dots,x_n]$  and let $I =
(G_1,\dots,G_{n+1})$ be the ideal of $n+1$ general 
forms of degrees $d_1 =\deg(G_1)\leq \dots \leq
d_{n+1}=\deg(G_{n+1})$. The Hilbert function  of $R/I$ is well
known (at least in characteristic zero -- see Remark \ref{characteristic
remark}), coming from a result of
R.\ Stanley
\cite{stanley} and of J.\ Watanabe \cite{watanabe} which implies that a
general Artinian complete intersection has the Strong Lefschetz Property,
and a very
    long-standing problem in
Commutative Algebra is to determine  the minimal free resolution
of $R/I$. In \cite{MMR2}, the first and second author
    gave
    the precise graded Betti numbers of $R/I$ in the following cases:
\label{MMR2 list}
\begin{itemize}
\item $n=3$.

\item $n=4$ and $\sum_{i=1}^5 d_i$ is even.

\item $n=4$, $\sum_{i=1}^{5} d_i$ is odd and $d_2 + d_3 + d_4 <
d_1 + d_5 + 4$.

\item $n$ is even and all generators have the same degree, $a$,
which is even.

\item $(\sum_{i=1}^{n+1} d_i) -n$ is even and $d_2 + \dots + d_n <
d_1 + d_{n+1} + n$.

\item $(\sum_{i=1}^{n+1} d_i) - n$ is odd, $n \geq 6 $ is even,
$d_2 + \dots+d_n < d_1 + d_{n+1} + n$ and $d_1 + \dots + d_n -
d_{n+1} - n \gg 0$.
\end{itemize}

As a nice application of Theorem \ref{res-rel-comp-gorcodim-n}, we
will enlarge the above list. Since the calculations are somewhat
complicated, we illustrate the method with an example before we
proceed to the general statement.

\begin{example} \label{exampleaci} Let $n=5$, $d_1=2$,
$d_2=d_3=d_4=4$, $d_5=5$ and $d_6=6$.  Consider   a general     Gorenstein Artinian graded algebra $R/G$ of embedding dimension
$5$, socle degree $8$ and relatively compressed with respect to a
complete intersection ideal $\fa =(G_1, G_2,G_3,G_4)$ with
$\deg(G_1)=2$, $\deg(G_2)=4$, $\deg(G_3)=4$ and $\deg(G_4)=4$.
    By Theorem
\ref{res-rel-comp-gorcodim-n} $A$ has a minimal free
    $R$-resolution of the following type:
    \[ 0 \rightarrow R(-13)
    \rightarrow R(-8)^{46}\oplus R(-11)\oplus R(-9)^3 \rightarrow R(-7)^{149}
    \rightarrow R(-6)^{149} \rightarrow
    \]
    \[R(-5)^{46}\oplus R(-2)\oplus R(-4)^3
     \rightarrow R \rightarrow R/G \rightarrow 0.
     \]

     \vskip 2mm
Hence, there exists a complete intersection $J\subset G$ with
generators of degrees 2,4,4,4,5. By a standard mapping cone
argument, the residual $I=[J:G]$ is an almost complete
intersection of type (2,4,4,4,5,6) and with the following minimal
free $R$-resolution \begin{equation} \label{example-aci} 0
\rightarrow R(-14)^{45} \rightarrow R(-13)^{146} \rightarrow
R(-12)^{150}\oplus R(-11)^3\oplus R(-10)^3 \rightarrow
\end{equation}
\[ R(-6)^3 \oplus R(-7)\oplus R(-8)^4 \oplus R(-9)^3
\oplus R(-10)^3 \oplus R(-11)^{46}\rightarrow \]
\[R(-2)\oplus
R(-4)^3\oplus R(-5)\oplus
    R(-6) \rightarrow R \rightarrow R/I \rightarrow 0.
\]

\vskip 2mm Since there are no non-Koszul ghost terms and the
graded Betti numbers are the smallest consistent with the Hilbert
function 1 5 14 30 52 75 92 95 79 45 0 of the general almost
complete intersection of type (2,4,4,4,5,6), the exact sequence
(\ref{example-aci}) gives us the minimal free $R$-resolution of
the general almost complete intersection of type (2,4,4,4,5,6).
\end{example}

The idea behind  Example \ref{exampleaci} leads to the following
result

\begin{theorem}\label{resol-aci} Let $I = (G_1,\dots,G_{n+1})$ be a
general almost
complete intersection in $R = k[x_1,\dots,x_n]$, with $d_i = \deg G_i$, $2
\leq d_1 \leq d_2 \dots \leq d_n \le d_{n+1} \leq (\sum_{i=1}^n
d_i) -n$, (the latter condition only assures that $I$ is not a
complete intersection,) and $\sum_{i=1}^{n+1} d_i -n$ even. Then,
$R/I$ has a minimal
free $R$-resolution of the form

\[ 0 \rightarrow F_1^\vee (-d) \rightarrow
\begin{array}{c}
K_1(\underline{\underline{d}})^\vee (-d) \\ \oplus \\ F_2^\vee
(-d)
\end{array}
\rightarrow
\begin{array}{c}
K_2(\underline{\underline{d}})^\vee (-d) \\ \oplus \\ F_3^\vee
(-d)
\end{array}
\rightarrow \cdots \hbox{\hskip 1.8in}
\]
\[\vbox{ \vskip .4in}
\hbox{\hskip 1.8in} \rightarrow
\begin{array}{c}
K_{n-2}(\underline{\underline{d}})^\vee (-d) \\ \oplus \\
F_{n-1}^\vee (-d)
\end{array}
\rightarrow \oplus _{i=1}^{n+1}R(-d_i) \rightarrow
R \rightarrow R/I \rightarrow 0
\]

\vskip 2mm \noindent  where $\underline{\underline{d}}=(d_1,
\cdots ,d_n)$, $d:=d_1+ \cdots +d_n$,
$K_i(\underline{\underline{d}})=K_{i}(d_1, \cdots ,d_n)$ and
$$F_i=K_{i}(\underline{d})^{\le (c/2)+i-1}\oplus R(-(c/2)-i)^{\alpha
_{i}(\underline{d},n,c/2)}\oplus (K_{n-i}(\underline{d})^{\le
(c/2)+n-i-1})^{\vee}(-c-n)$$

\vskip 2mm \noindent with $c=(\sum _{i=1}^{n} d_i)-d_{n+1}-n$,
$r=min(n,max\{ i / d_i\le c/2 \})$ $\underline{d}=(d_1, \cdots
,d_r)$, $K_i(\underline{d})=K_{i}(d_1, \cdots ,d_r)$, $\alpha
_{i}(\underline{d},n,c/2)=\alpha _{n-i}(\underline{d},n,c/2)$ and
$\alpha _{i}(\underline{d},n,c/2)$  determined by the Hilbert
function of $R/(G_1,\cdots ,G_r)$.
\end{theorem}

\begin{proof}
Consider   a general 
    Gorenstein Artinian graded algebra $R/G$ of embedding dimension
$n$, socle degree $c=(\sum _{i=1}^{n} d_i)-d_{n+1}-n$ and
relatively compressed with respect to a complete intersection
ideal $\fa =(G_1, \cdots ,G_r)$ with $\deg(G_i)=d_i$ and
$r=\min(n,\max\{ i | d_i\le c/2 \})$. So,
$h_{R/G}(t)=\min\{h_{R/(G_1,\cdots , G_r)}(t),h_{R/(G_1,\cdots ,
G_r)}(c-t)\}$.

Since by hypothesis the socle degree $c$ of $R/G$ is even, we can
apply
  Theorem
\ref{res-rel-comp-gorcodim-n} and  we get that $R/G$ has a minimal
free
    $R$-resolution of the following type:
    \[
0 \rightarrow R(-c-n) \rightarrow F_{n-1} \rightarrow \cdots
\rightarrow F_2 \rightarrow F_1 \rightarrow R \rightarrow R/G
\rightarrow 0
\]

\noindent where  for all $i=1, \cdots , n-1$, we have

$$F_i=K_{i}(\underline{d})^{\le (c/2)+i-1}\oplus R(-(c/2)-i)^{\alpha
_{i}(\underline{d},n,c/2)}\oplus (K_{n-i}(\underline{d})^{\le
(c/2)+n-i-1})^{\vee}(-c-n)$$

\vskip 2mm \noindent with  $\alpha
_{i}(\underline{d},n,c/2)=\alpha _{n-i}(\underline{d},n,c/2)$ and
$\alpha _{i}(\underline{d},n,c/2)$ is determined by the Hilbert
function of $R/G$.

    \vskip 2mm Hence, there exists a complete
intersection $J\subset G$ with generators of degrees $d_1, \cdots
, d_n$.  The minimal free $R$-resolution of $R/J$ is  given by the
Koszul resolution:
\[
0 \rightarrow K_n (\underline{\underline{d}})\rightarrow
K_{n-1}(\underline{\underline{d}}) \rightarrow \cdots \rightarrow
K_2 (\underline{\underline{d}})\rightarrow K_1
(\underline{\underline{d}})\rightarrow R \rightarrow R/J
\rightarrow 0.
\]

By a standard mapping cone argument applied to the diagram
\[
\begin{array}{ccccccccccccccccccccc}
0 & \rightarrow & R(-d) & \rightarrow &
K_{n-1}(\underline{\underline{d}}) & \rightarrow & \cdots &
\rightarrow & K_2 (\underline{\underline{d}})& \rightarrow &
K_1(\underline{\underline{d}}) & \rightarrow & R & \rightarrow & R/J
& \rightarrow & 0 \\
&& \downarrow && \downarrow &&&& \downarrow && \downarrow &&
\downarrow && \downarrow \\ 0 & \rightarrow & R(-c-n) &
\rightarrow & F_{n-1} & \rightarrow & \cdots & \rightarrow & F_2 &
\rightarrow & F_1 & \rightarrow & R & \rightarrow & R/G &
\rightarrow & 0
\end{array}
\]
the residual $I=[J:G]$ is an almost complete intersection of type
$$(d_1, \cdots , d_n,-c-n+\sum _{i=1}^{n} d_i )=(d_1, \cdots , d_n,
d_{n+1})$$ and with the following minimal free $R$-resolution

\begin{equation}\label{res-aci}
0 \rightarrow F_1^\vee (-d) \rightarrow
\begin{array}{c}
K_1(\underline{\underline{d}})^\vee (-d) \\ \oplus \\ F_2^\vee
(-d)
\end{array}
\rightarrow
\begin{array}{c}
K_2(\underline{\underline{d}})^\vee (-d) \\ \oplus \\ F_3^\vee
(-d)
\end{array}
\rightarrow \cdots \rightarrow \oplus _{i=1}^{n+1}R(-d_i)
\rightarrow R \rightarrow R/I \rightarrow 0.
\end{equation}

\vskip 2mm Since there are no non-Koszul ghost terms and the
graded Betti numbers are the smallest consistent with the Hilbert
function
\[
h_{R/I'}(\ell ) = \left [ h_{R/J} \left ( (\sum_{i=1}^n d_i) - n - \ell
\right ) - h_{R/J}
\left ( (\sum_{i=1}^n d_i) - n - \ell - d_{n+1} \right ) \right ]_+ \]
    (where $[x]_+$ denotes the maximum of $x$ and 0) of the general
almost complete
intersection ideal $I'\subset k[x_1, \cdots ,x_n] $ of type $(d_1,
\cdots , d_n, d_{n+1})$, the exact sequence (\ref{res-aci}) gives
us  the minimal free $R$-resolution of the general almost complete
intersection of type $(d_1, \cdots , d_n, d_{n+1})$.
\end{proof}

\vskip 4mm

\begin{remark}
We point out that there are many new cases covered by Theorem
\ref{resol-aci}, that were not known previously (and in particular not in
the list on page \pageref{MMR2 list}).
The most natural remaining open question is to determine the minimal free
resolution of an ideal of $n+1$ general forms of degree $a$ in $n$
variables, when either $n$ is odd or $n$ is even and $a$ is odd.
\end{remark}

\vskip 2mm

\end{document}